\documentclass[a4paper,12pt]{article}
\usepackage[utf8]{inputenc} 
\usepackage[english]{babel}
\usepackage{amsfonts,amssymb}
\title{Vector Lyapunov-Krasovsky functionals, higher derivatives of Lyapunov-Krasovsky functionals and boundedness of solutions of functional-differential equations of retarded type\footnote{This research was supported by the Russian Science Foundation (project No. 25-21-00830)}}
\author{K.S. Lapin, S.V. Lapin}
\date{}
\begin{document}
\maketitle

\begin{abstract}
In present paper, for functional-differential equations of retarded type, the method of vector Lyapunov-Krasovsky functionals and the method of higher derivatives of Lyapunov-Krasovsky functionals are developed. On the basis of this methods, sufficient conditions for various types of uniform boundedness of solutions and sufficient conditions for various types of ultimate-uniform boundedness of solutions are obtained.

\end{abstract}
\vspace{0.5cm}

\centerline{\bf 1. Introduction}
\vspace{0.5cm}

In the work of T. Yoshizawa \cite{Io2}, on the basis of the method of Lyapunov functions \cite{L}, a theory of boundedness of solutions of systems of ordinary differential equations (ODEs) was created. In order to emphasize that this theory is constructed using Lyapunov functions, in what follows this theory will be called the LF-theory of Yoshizawa boundedness (LF -- Lyapunov functions). A generalization of the LF-theory of Yoshizawa boundedness on the basis of the method of vector Lyapunov functions \cite{A}, hereinafter called the VLF-theory of Yoshizawa boundedness (VLF -- vector Lyapunov functions), was carried out in the monograph of V.M. Matrosov \cite{M}. Further, in the works of K.S. Lapin (see, for example, \cite{Lap1}-\cite{Lap3}), on the basis of the methods of Lyapunov functions and vector Lyapunov functions, the VLF-theory of partial boundedness of
solutions of systems of ODEs with partial control of initial conditions was developed. This theory generalizes the above-mentioned VLF-theory of Yoshizawa boundedness of solutions of systems of ODEs in the sense that if the part of the coordinates of solutions under study coincides with the complete set of coordinates of solutions and, moreover, the controlled part of the coordinates of the initial conditions of solutions coincides with the complete set of coordinates of the initial conditions of solutions, then the VLF-theory of boundedness developed in the works of K.S. Lapin becomes the VLF-theory of Yoshizawa boundedness. On the other hand, in the work of T. Yoshizawa \cite{Io1}, on the basis of the method of Lyapunov-Krasovsky functionals \cite{Kra}, a theory of boundedness of solutions of functional-differential equations (FDEs) of retarded type was created. In what follows this theory will be called the LKF-theory of Yoshizawa boundedness for FDEs (LKF -- Lyapunov-Krasovsky functionals). This LKF-theory of Yoshizawa boundedness for FDEs generalizes the LF-theory of Yoshizawa boundedness for ODEs in the sense that in this theory, instead of a finite-dimensional space of initial coordinate positions of solutions, an infinite-dimensional functional space of prehistories of solutions is used. The main difficulties in the development of the LKF-theory of Yoshizawa boundedness for FDEs, as was shown in \cite{Io1} (see also \cite{H}), are connected precisely with the infinite-dimensionality of the functional space of prehistories of solutions. In light of what has been said above, an interesting and important problem for theory and applications naturally arises: to develop the method of vector Lyapunov-Krasovsky functionals and, on the basis of this method, to construct the VLKF-theory of partial boundedness of solutions of FDEs of retarded type with partial control of prehistories of solutions (VLKF -- vector Lyapunov-Krasovsky functionals). A natural requirement for the construction of such a VLKF-theory is that this theory, on the one hand, generalizes the LKF-theory of Yoshizawa for FDEs, and, on the other hand, generalizes the VLF-theory of K.S. Lapin of partial boundedness of solutions of systems of ODEs with partial control of initial conditions.
\vspace{0.5cm}

\centerline{\bf 2. Vector Lyapunov-Krasovsky functionals and}
\centerline{\bf higher derivatives of Lyapunov-Krasovsky functionals}
\vspace{0.5cm}

Let $C([-h;0],\mathbb{R}^n)$ be the linear space of all continuous mappings of the form $\varphi:[-h;0]\to\mathbb{R}^n$, $n\geqslant 1$, considered with the norm
$\|\varphi\|=\max_{-h\leqslant\vartheta\leqslant 0}\|\varphi(\vartheta)\|$,
where $h>0$ is a fixed number and $\|\varphi(\vartheta)\|$ is the Euclidean norm of the vector $\varphi(\vartheta)\in\mathbb{R}^n$.
It is easy to see that any continuous mapping $x:[-h;A)\to\mathbb{R}^n$, $x(t)\in\mathbb{R}^n$, $-h\leqslant t< A$, where $A>0$ or $A=+\infty$, determines a continuous mapping $\overline{x}:[0;A)\to C([-h;0],\mathbb{R}^n)$, $\overline{x}(t)=x_t\in C([-h;0],\mathbb{R}^n)$, $0\leqslant t<A$, where the mapping $x_t:[-h;0]\to\mathbb{R}^n$ is defined by the formula
$x_t(\vartheta)=x(t+\vartheta)$, $-h\leqslant \vartheta\leqslant 0$.
Obviously, $x_t(0)=x(t)$ and $x_t(-h)=x(t-h)$.

Consider, following \cite{Io1} (see also \cite{H}), the functional-differential equation of retarded type
$$\dot{x}(t)=F(t,x_t),\quad t\in\mathbb{R}^+\equiv [0;+\infty),\eqno(2.1)$$
where $x(t)\in\mathbb{R}^n$, $-h\leqslant t<+\infty$, $F:\mathbb{R}^+\times C([-h;0],\mathbb{R}^n)\to \mathbb{R}^n$ is any continuous mapping, $\dot{x}(t)$ is the right-hand derivative of the mapping $x(t)$ and the family of mappings $x_t:[-h;0]\to \mathbb{R}^n$, $t\in\mathbb{R}^+$, is determined in the above-mentioned way from the mapping $x(t)$, $-h\leqslant t<+\infty$. The solution $x(t)$ of equation $(2.1)$ with initial conditions $(t_0,\varphi_0)\in \mathbb{R}^+\times C([-h;0],\mathbb{R}^n)$ will in what follows be denoted by $x(t,t_0,\varphi_0)$. The mapping $x_t$ corresponding to the solution $x(t,t_0,\varphi_0)$ will be denoted by $x_t(t_0,\varphi_0)$. The mapping $\varphi_0:[-h;0]\to\mathbb{R}^n$ will in what follows be called the prehistory of the solution $x(t,t_0,\varphi_0)$ from the time moment $t_0-h$ to the time moment $t_0$ or, more briefly, the prehistory of this solution. Moreover, in what follows we shall assume that in equation $(2.1)$ the continuous mapping $F(t,\varphi)$ satisfies a local Lipschitz condition with respect to the variable $\varphi$.

Examples of functional-differential equations of the form $(2.1)$ are given by differential equations with positive constant time delay. Indeed, consider the differential equation
$$\dot{x}(t)=f(t,x(t),x(t-h))\eqno(2.2)$$
with constant delay $h>0$, where $x(t)\in\mathbb{R}^n$, $t_0-h\leqslant t<+\infty$, $t_0\in\mathbb{R}^+$, $\dot{x}(t)$ is the right-hand derivative of the mapping $x(t)$, $t\geqslant t_0$, and the mapping $f:\mathbb{R}^+\times\mathbb{R}^n\times\mathbb{R}^n\to \mathbb{R}^n$ is any continuous mapping satisfying a local Lipschitz condition with respect to the variables $(x_1,x_2)\in\mathbb{R}^n\times\mathbb{R}^n$. For $n\geqslant 2$, equations of the form $(2.2)$ are also called systems of differential equations with delay. It is easy to see that every differential equation $(2.2)$ with constant delay $h>0$ determines a functional-differential equation of retarded type $\dot{x}(t)=F(t,x_t)$,
where the continuous mapping $F:\mathbb{R}^+\times C([-h;0],\mathbb{R}^n)\to \mathbb{R}^n$ is defined by the formula
$F(t,\varphi)=f(t,\varphi(0),\varphi(-h))$.

Let any continuous functional $U:\mathbb{R}^+\times C([-h;0],\mathbb{R}^n)\to\mathbb{R}$ be given. Consider for this functional
and the solution
$x(\tau,t,\varphi)$ of equation $(2.1)$ with any initial conditions $(t,\varphi)$ the continuous function $u(\tau)=U(\tau,x_\tau(t,\varphi))$, $\tau\geqslant t$. The derivative of this functional $U(t,\varphi)$
along equation $(2.1)$ is defined as the mapping $\dot{U}_{(2.1)}:\mathbb{R}^+\times C([-h;0],\mathbb{R}^n)\to\mathbb{R}\cup\{\pm\infty\}$ given by the formula
$$\dot{U}_{(2.1)}(t,\varphi)=u'^+(\tau)|_{\tau=t}=u'^+(t)=
\overline{\lim_{\delta\to\, 0^+}}\,\,\frac{1}{\delta}\left\{U(t+\delta,x_{t+\delta}(t,\varphi))-U(t,\varphi)\right\},$$
where $u'^+(\tau)$ is the upper right derivative (see, for example, \cite{R}) of the continuous function $u(\tau)$.

Let us now proceed to the consideration of vector functional and their derivatives along equation $(2.1)$

Let an arbitrary continuous vector functional
$$V:\mathbb{R}^+\times C([-h;0],\mathbb{R}^n)\to\mathbb{R}^l,\quad l\geqslant 1,$$
$$V(t,\varphi)=(V_1(t,\varphi),\dots,V_l(t,\varphi)),\quad (t,\varphi)\in\mathbb{R}^+\times C([-h;0],\mathbb{R}^n)$$ be given.
The derivative $\dot{V}_{(2.1)}(t,\varphi)$ of the vector functional $V(t,\varphi)$ along
equation $(2.1)$ is defined by the equality
$\dot{V}_{(2.1)}(t,\varphi)=((\dot{V}_{(2.1)})_1(t,\varphi),\dots,(\dot{V}_{(2.1)})_l(t,\varphi))$, where
$(\dot{V}_{(2.1)})_i(t,\varphi)=(\dot{V}_i)_{(2.1)}(t,\varphi)$ is the derivative of the functional
$V_i(t,\varphi)$, $1\leqslant i\leqslant l$, along equation $(2.1)$.

Let now a continuous vector function
$$f:\mathbb{R}^+\times\mathbb{R}^l\to\mathbb{R}^l,\quad l\geqslant 1,\quad f(t,\xi)=(f_1(t,\xi),\dots,f_l(t,\xi)),\quad (t,\xi)\in\mathbb{R}^+\times\mathbb{R}^l$$ be given.
In what follows we shall use the
notation $f(t,\xi)\in W$ if $f(t,\xi)$ satisfies the
Wazewski condition \cite{W}, which consists in the fact that for any fixed numbers $1\leqslant
s\leqslant l$ and $t\in\mathbb{R}^+$ the function $f_s(t,\xi)$ is non-decreasing in the variables
$\xi_1,\dots,\xi_{s-1},\xi_{s+1},\dots,\xi_l$, i.e., from the conditions
$\xi_i\leqslant\eta_i$, $1\leqslant i\leqslant l$, $i\not=s$,
and $\xi_s=\eta_s$ it follows that $f_s(t,\xi)\leqslant f_s(t,\eta)$.
Note
that for $l=1$ the condition $f(t,\xi)\in W$ degenerates, therefore for any continuous
function $f(t,\xi)$, where $(t,\xi)\in\mathbb{R}^+\times\mathbb{R}$, we shall agree
to assume that $f(t,\xi)\in W$.

For any vectors
$\xi=(\xi_1,\dots,\xi_l)\in\mathbb{R}^l$ and
$\eta=(\eta_1,\dots,\eta_l)\in\mathbb{R}^l$ in what follows the
notation $\xi\leqslant\eta$ is used if $\xi_i\leqslant\eta_i$ for any
$1\leqslant i\leqslant l$.

{\bf Definition 2.1}. A continuous vector functional
$$V:\mathbb{R}^+\times C([-h;0],\mathbb{R}^n)\to\mathbb{R}^l,\quad l\geqslant 1,$$
$$V(t,\varphi)=(V_1(t,\varphi),\dots,V_l(t,\varphi)),\quad (t,\varphi)\in\mathbb{R}^+\times C([-h;0],\mathbb{R}^n),$$
for which the derivative $\dot{V}_{(2.1)}(t,\varphi)$ along equation $(2.1)$ is also a continuous vector functional, will in what follows be called a vector Lyapunov-Krasovsky functional for equation $(2.1)$ if there exists a continuous vector function
$$f:\mathbb{R}^+\times\mathbb{R}^l\to\mathbb{R}^l,\quad l\geqslant 1,$$
$$f(t,\xi)=(f_1(t,\xi),\dots,f_l(t,\xi)),\quad (t,\xi)\in\mathbb{R}^+\times\mathbb{R}^l,\quad f(t,\xi)\in W,$$
for which for all $(t,\varphi)\in\mathbb{R}^+\times C([-h;0],\mathbb{R}^n)$ the inequality
$$\dot{V}_{(2.1)}(t,\varphi)\leqslant f(t,V(t,\varphi)).\eqno(2.3)$$ is valid.

{\bf Definition 2.2}. The comparison system for equation $(2.1)$ and the vector Lyapunov-Krasovsky functional $V(t,\varphi)$ for this equation, which satisfies inequality $(2.3)$, will be called the system of ordinary differential equations $$\dot{\xi}=f(t,\xi),\quad\xi\in\mathbb{R}^l,\quad f(t,\xi)\in W,\eqno(2.4)$$
where $f(t,\xi)$ is the vector function which stands in the right-hand side of inequality $(2.3)$. In what follows, if this does not cause misunderstandings, system $(2.4)$ will, for brevity, be called the comparison system for equation $(2.1)$.

Since the right-hand
side of the comparison system $(2.4)$ for equation $(2.1)$ is a continuous vector function, the uniqueness of the solution of the Cauchy problem for this system may fail. However, as was shown by T. Wazewski \cite{W}, for any system of ordinary differential equations $\dot{\xi}=g(t,\xi)$, $\xi\in\mathbb{R}^l$, $l\geqslant 1$, where $g(t,\xi)\in W$,
among all solutions of this system with any given initial conditions $(t_0,\xi_0)$, one can always choose the upper
solution $\overline{\xi}(t,t_0,\xi_0)$, i.e., the solution for which for all $t\geqslant t_0$
the inequality
$\xi(t,t_0,\xi_0)\leqslant\overline{\xi}(t,t_0,\xi_0)$ is valid, where $\xi(t,t_0,\xi_0)$ is any solution of the system $\dot{\xi}=g(t,\xi)$ with initial conditions $(t_0,\xi_0)$. From this result of Wazewski it follows that the comparison system $(2.4)$ for equation $(2.1)$ always has an upper solution $\overline{\xi}(t,t_0,\xi_0)$ with arbitrary initial conditions $(t_0,\xi_0)\in\mathbb{R}^+\times\mathbb{R}^l$.

Let us now establish the connection between the solutions of equation $(2.1)$ and the upper solutions of the comparison system $(2.4)$ for this equation.

{\bf Theorem 2.1}. For any vector Lyapunov-Krasovsky functional $V(t,\varphi)$ with comparison system $(2.4)$ for equation $(2.1)$ and any solution $x(t,t_0,\varphi_0)$ of equation $(2.1)$ with any initial conditions $(t_0,\varphi_0)\in \mathbb{R}^+\times C([-h;0],\mathbb{R}^n)$, for all $t\geqslant t_0$ the inequality
$$V(t,x_t(t_0,\varphi_0))\leqslant\overline{\xi}(t,t_0,V(t_0,\varphi_0)),\eqno(2.5)$$
is valid, where $\overline{\xi}(t,t_0,V(t_0,\varphi_0))$ is the upper solution of the comparison system $(2.4)$ for equation $(2.1)$ with initial conditions $(t_0,V(t_0,\varphi_0))\in\mathbb{R}^+\times\mathbb{R}^l$.

{\bf Proof}. Let us first recall that in the work of T. Wazewski \cite{W} it was shown that if some continuously differentiable vector function $u:\mathbb{R}^+\to\mathbb{R}^l$, $l\geqslant 1$,
satisfies for all $t\in\mathbb{R}^+$ the inequality
$u'(t)\leqslant g(t,u(t))$, where
$g:\mathbb{R}^+\times\mathbb{R}^l\to\mathbb{R}^l$, $l\geqslant 1$, is some continuous vector function
satisfying the condition $g(t,\xi)\in W$, then
for the upper solution $\overline{\xi}(t,t_0,u(t_0))$ of the equation $\dot{\xi}=g(t,\xi)$ with initial conditions $(t_0,u(t_0))\in\mathbb{R}^+\times\mathbb{R}^l$ the inequality $u(t)\leqslant\overline{\xi}(t,t_0,u(t_0))$ is valid for all $t\geqslant t_0$.
Let us now connect this result of Wazewski with vector Lyapunov-Krasovsky functionals. Let any vector Lyapunov-Krasovsky functional $V(t,\varphi)$ with comparison system $(2.4)$ for equation $(2.1)$ be given. Let us take as the above-mentioned vector function $u(t)$ the vector function $v(t)=V(t,x_t(t_0,\varphi_0))$, $t\in\mathbb{R}^+$, where $x(t,t_0,\varphi_0)$ is the solution of equation $(2.1)$ with any initial conditions $(t_0,\varphi_0)\in\mathbb{R}^+\times
C([-h;0],\mathbb{R}^n)$. Since the derivative $\dot{V}_{(2.1)}(t,\varphi)$ of the vector Lyapunov-Krasovsky functional $V(t,\varphi)$ along equation $(2.1)$ is by definition a continuous vector functional, the vector function $v(t)$ is continuously differentiable. From inequality $(2.3)$ it follows that for any point $t_0\in\mathbb{R}^+$ the inequality $v'(t_0)\leqslant f(t_0,v(t_0))$ is valid, since $v'(t_0)=\dot{V}(t_0,\varphi_0)$ and $v(t_0)=V(t_0,\varphi_0)$. Since the point $t_0\in\mathbb{R}^+$ is arbitrary, for all $t\in\mathbb{R}^+$ the inequality $v'(t)\leqslant f(t,v(t))$ is valid, where $f(t,\xi)\in W$. Applying in this situation the above-mentioned result of Wazewski, we obtain the inequality $v(t)\leqslant\overline{\xi}(t,t_0,v(t_0))$ valid for all $t\geqslant t_0$, where $\overline{\xi}(t,t_0,v(t_0))$ is the upper solution of the comparison system $(2.4)$ for equation $(2.1)$ with initial conditions $(t_0,v(t_0))\in\mathbb{R}^+\times\mathbb{R}^l$. From this inequality, taking into account that $v(t)=V(t,x_t(t_0,\varphi_0))$ and $v(t_0)=V(t_0,\varphi_0)$, we obtain the required inequality $V(t,x_t(t_0,\varphi_0))\leqslant\overline{\xi}(t,t_0,V(t_0,\varphi_0))$,
where $\overline{\xi}(t,t_0,V(t_0,\varphi_0))$ is the upper solution of the comparison system $(2.4)$ for equation $(2.1)$ with initial conditions $(t_0,V(t_0,\varphi_0))\in\mathbb{R}^+\times\mathbb{R}^l$.~~~$\blacksquare$

Let us now proceed to the consideration of higher derivatives of continuous functionals along equation $(2.1)$.

Let an arbitrary continuous functional
$U:\mathbb{R}^+\times C([-h;0],\mathbb{R}^n)\to\mathbb{R}$ be given. The derivative $\dot{U}_{(2.1)}(t,\varphi)$ of this functional along equation $(2.1)$ will in what follows be denoted by $U^{(1)}_{(2.1)}(t,\varphi)$ and called the first-order derivative of the functional $U(t,\varphi)$ along system $(2.1)$. If the first-order derivative $U^{(1)}_{(2.1)}(t,\varphi)$ of the functional $U(t,\varphi)$ along equation $(2.1)$ is a continuous functional
$U^{(1)}_{(2.1)}:\mathbb{R}^+\times C([-h;0],\mathbb{R}^n)\to\mathbb{R}$, then the derivative of this functional along equation $(2.1)$ will be denoted by $U^{(2)}_{(2.1)}(t,\varphi)$ and called the second-order derivative of the functional $U(t,\varphi)$ along system $(2.1)$. Analogously, the $s$-th order derivative $U^{(s)}_{(2.1)}(t,\varphi)$, $s\geqslant 1$, is defined by the following recurrence formula: $$U^{(s)}_{(2.1)}(t,\varphi)=(U^{(s-1)}_{(2.1)})^{(1)}(t,\varphi),\quad U^{(0)}_{(2.1)}(t,\varphi)=U(t,\varphi),\quad U^{(1)}_{(2.1)}(t,\varphi)=\dot{U}_{(2.1)}(t,\varphi).$$
For $s>1$, the derivatives $U^{(s)}_{(2.1)}(t,\varphi)$ are usually called the higher derivatives of the functional $U(t,\varphi)$ along equation $(2.1)$.

Let now a function $g(t,\sigma)$ be given, where $(t,\sigma)=(t,\sigma_1,\dots,\sigma_l)\in\mathbb{R}^+\times\mathbb{R}^l$, $l\geqslant 1$.
In what follows we shall say that the function $g(t,\sigma)$ is non-decreasing in the variables $\sigma_1,\dots,\sigma_{l-1}$
if from the conditions $\sigma_s\leqslant\eta_s$, $1\leqslant s\leqslant
l-1$, and $\sigma_l=\eta_l$ it follows that $g(t,\sigma)\leqslant
g(t,\eta)$.

{\bf Definition 2.3}. A continuous functional
$U:\mathbb{R}^+\times C([-h;0],\mathbb{R}^n)\to\mathbb{R}$
for which the $l$-th order derivative $U^{(l)}_{(2.1)}(t,\varphi)$ along equation $(2.1)$ is a continuous vector functional $U^{(l)}_{(2.1)}:\mathbb{R}^+\times C([-h;0],\mathbb{R}^n)\to\mathbb{R}$, will be called a Lyapunov-Krasovsky functional with an $l$-th order derivative along equation $(2.1)$ if there exists a continuous function
$g(t,\sigma)$, where $(t,\sigma)=(t,\sigma_1,\dots,\sigma_l)\in\mathbb{R}^+\times\mathbb{R}^l$, non-decreasing in the variables $\sigma_1,\dots,\sigma_{l-1}$,
for which for all $(t,\varphi)\in\mathbb{R}^+\times C([-h;0],\mathbb{R}^n)$ the inequality
$$U^{(l)}_{(2.1)}(t,\varphi)\leqslant g(t,U(t,\varphi),\dot{U}_{(2.1)}(t,\varphi),\dots,U^{(l-1)}_{(2.1)}(t,\varphi)).\eqno(2.6)$$ is valid.

{\bf Definition 2.4}. The $l$-th order comparison equation for equation $(2.1)$ and the Lyapunov-Krasovsky functional $U(t,\varphi)$ with an $l$-th order derivative along equation $(2.1)$, which satisfies inequality $(2.6)$, will be called the ordinary differential equation of $l$-th order
$$\psi^{(l)}=g(t,\psi,\dot{\psi},\dots,\psi^{(l-1)}),\quad\psi\in\mathbb{R},\eqno(2.7)$$
where $g(t,\sigma)$ is the continuous function, non-decreasing in the variables $\sigma_1,\dots,\sigma_{l-1}$, which stands in the right-hand side of inequality $(2.6)$. In what follows, if this does not cause misunderstandings, equation $(2.7)$ will, for brevity, be called the $l$-th order comparison equation for equation $(2.1)$.

Consider for the ordinary differential equation of $l$-th order $(2.7)$ the corresponding system of ordinary differential equations of first order
$$\left.\begin{array}{c}
\dot{\xi}=f(t,\xi),\quad (t,\xi)=(t,\xi_1,\dots,\xi_l)\in\mathbb{R}^+\times\mathbb{R}^l,\\
f(t,\xi)=(f_1(t,\xi),\dots,f_l(t,\xi)),\quad\xi_1=\psi,\quad \dot{\xi}_1=f_1(t,\xi)=\xi_2,\\$$
\dots,\quad\dot{\xi}_{l-1}=f_{l-1}(t,\xi)=\xi_l,\quad \dot{\xi}_l=f_l(t,\xi)=g(t,\xi_1,\dots,\xi_l).\\
\end{array}\right\}\eqno(2.8)$$

Let us now establish the connection between higher derivatives of Lyapunov-Krasovsky functionals and vector Lyapunov-Krasovsky functionals.

{\bf Theorem 2.2}. Let any Lyapunov-Krasovsky functional $U(t,\varphi)$ with an $l$-th order derivative along equation $(2.1)$ and with comparison equation $(2.7)$ for equation $(2.1)$ be given. Then the vector functional
$$V(t,\varphi)=(V_1(t,\varphi),\dots,V_l(t,\varphi)),\quad
(t,\varphi)\in\mathbb{R}^+\times C([-h;0],\mathbb{R}^n),$$ defined by the formulas
$$V_1(t,\varphi)=U(t,\varphi),\quad V_2(t,\varphi)=\dot{U}_{(2.1)}(t,\varphi),\quad\dots,\quad V_l(t,\varphi)
=U^{(l-1)}_{(2.1)}(t,\varphi),\eqno(2.9)$$
is a vector Lyapunov-Krasovsky functional for equation $(2.1)$ with comparison system $(2.8)$ for equation $(2.1)$.

{\bf Proof}. Let us first note that the right-hand side
$f(t,\xi)$ of system $(2.8)$ satisfies the condition $f(t,\xi)\in W$, i.e.,
for every $1\leqslant i\leqslant l$ the continuous function $f_i(t,\xi)$
is non-decreasing in the variables $\xi_1,\dots,\xi_{i-1},\xi_{i+1},\dots,\xi_l$. Indeed,
for $1\leqslant i\leqslant l-1$ this is obvious because $f_i(t,\xi)=\xi_{i+1}$. For $i=l$ the continuous function
$f_l(t,\xi)=g(t,\xi_1,\dots,\xi_l)$ is, by the definition of the comparison equation $(2.7)$, non-decreasing in the variables
$\xi_1,\dots,\xi_{l-1}$. Let us now show that for the vector functional $V(t,\varphi)=(V_1(t,\varphi),\dots,V_l(t,\varphi))$,
defined by formulas $(2.9)$, and the continuous vector function $f(t,\xi)=(f_1(t,\xi),\dots,f_l(t,\xi))$, defined by formulas $(2.8)$, the inequality $(2.3)$ holds, i.e.,
for every $1\leqslant i\leqslant l$ and all $(t,\varphi)\in\mathbb{R}^+\times C([-h;0],\mathbb{R}^n)$ the inequality $(\dot{V}_i)_{(2.1)}(t,\varphi)\leqslant f_i(t,V(t,\varphi))$ is valid.
Indeed, for $1\leqslant i\leqslant
l-1$ we have the equality$(\dot{V}_i)_{(2.1)}(t,\varphi)=f_i(t,V(t,\varphi))$ for all $(t,\varphi)\in\mathbb{R}^+\times C([-h;0],\mathbb{R}^n)$ because
for all $(t,\varphi)\in\mathbb{R}^+\times C([-h;0],\mathbb{R}^n)$ the equalities $(\dot{V}_i)_{(2.1)}(t,\varphi)=U^{(i)}_{(2.1)}(t,\varphi)$ and $f_i(t,V(t,\varphi))=V_{i+1}(t,\varphi)=U^{(i)}_{(2.1)}(t,\varphi)$ are valid.
For $i=l$ we have the inequality $(\dot{V}_l)_{(2.1)}(t,\varphi)\leqslant f_l(t,V(t,\varphi))$ for all $(t,\varphi)\in\mathbb{R}^+\times C([-h;0],\mathbb{R}^n)$, since $(\dot{V}_l)_{(2.1)}(t,\varphi)=U^{(l)}_{(2.1)}(t,\varphi)$ and
$f_l(t,V(t,x))=g(t,V(t,\varphi))$ and, moreover, by the condition of the theorem being proved inequality $(2.6)$ is valid. Thus, it is shown that the vector functional $V(t,\varphi)$, defined by formulas $(2.9)$, is a vector Lyapunov-Krasovsky functional for equation $(2.1)$ with comparison system $(2.8)$ for equation $(2.1)$.~~~$\blacksquare$
\vspace{0.5cm}

\centerline{\bf 3. Vector Lyapunov-Krasovsky functionals, higher derivatives}
\centerline{\bf of Lyapunov-Krasovsky functionals and uniform boundedness}
\centerline{\bf of solutions of functional-differential equations of retarded type}
\vspace{0.5cm}

Let us first recall \cite{Io1} that solutions of equation $(2.1)$ are said to be uniformly bounded if for every number
$\alpha\geqslant 0$ there exists a number $\beta=\beta(\alpha)>\alpha$ such that for every solution $x(t,t_0,\varphi_0)$, $\|\varphi_0\|_C\leqslant\alpha$,
of equation $(2.1)$ the condition
$\|x_t(t_0,\varphi_0)\|_C\leqslant\beta$ is satisfied for all $t\geqslant t_0$.

It is easy to see that the concept of uniform boundedness of solutions of equation $(2.1)$ is equivalent to the fact that for every number
$\alpha\geqslant 0$ there exists a $\beta=\beta(\alpha)>\alpha$ such that for every solution $x(t,t_0,\varphi_0)$,
$\|\varphi_0\|_C\leqslant\alpha$, of equation $(2.1)$ the condition $\|x(t,t_0,\varphi_0)\|\leqslant\beta$ is satisfied for all $t\geqslant t_0$.

In what follows, for every element $x=(x_1,\dots,x_n)\in\mathbb{R}^n$, $n\geqslant 1$, and any
integer $1\leqslant k\leqslant n$
we shall use the notation $y=(x_1,\dots,x_k)\in\mathbb{R}^k$. Moreover, in what follows, for every $\varphi=(\varphi_1,\dots,\varphi_n)\in C([-h;0],\mathbb{R}^n)$ and any integer $1\leqslant s\leqslant n$ we shall use the notation $\varphi^{(s)}=(\varphi_1,\dots,\varphi_s)\in C([-h;0],\mathbb{R}^s)$.

Let us introduce the concept of uniform partial boundedness of solutions of equation $(2.1)$, which generalizes the above-mentioned concept of uniform boundedness of solutions of equation $(2.1)$.

{\bf Definition 3.1}. Solutions of equation $(2.1)$ are said to be uniformly bounded with respect to a part of the variables
$y=(x_1,\dots,x_k)$, $1\leqslant k\leqslant n$, or, more briefly, uniformly
$y$-bounded, if for every number
$\alpha\geqslant 0$ there exists a number $\beta=\beta(\alpha)>\alpha$ such that for
every solution $x(t,t_0,\varphi_0)$, $\|\varphi_0\|_C\leqslant \alpha$, of equation $(2.1)$ the condition
$\|y_t(t_0,\varphi_0)\|_C\leqslant\beta$ is satisfied for all $t\geqslant t_0$, where $y_t(t_0,\varphi_0)=x_t(t_0,\varphi_0)^{(k)}$.

It is clear that if in Definition 3.1 we set $k=n$, then
we obtain the definition of uniform boundedness of solutions of equation $(2.1)$ considered above.

It is easy to see that the concept of uniform $y$-boundedness of solutions of equation $(2.1)$ is equivalent to the fact that for every number $\alpha\geqslant 0$ there exists a number $\beta=\beta(\alpha)>\alpha$ such that for every solution $x(t,t_0,\varphi_0)$
of equation $(2.1)$, where $\|\varphi_0\|_C\leqslant \alpha$, the condition $\|y(t,t_0,\varphi_0)\|\leqslant\beta$ is satisfied for all $t\geqslant t_0$.

Let us introduce the concept of uniform partial boundedness of solutions of equation $(2.1)$ with partial control of prehistories, which generalizes the concept of uniform $y$-boundedness of solutions of equation $(2.1)$ introduced above.

{\bf Definition 3.2}. Solutions of equation $(2.1)$ are said to be uniformly bounded with respect to a part of the variables
$y=(x_1,\dots,x_k)$, $1\leqslant k\leqslant n$, with partial control of prehistories $\varphi^{(m)}_0=((\varphi_0)_1,\dots,(\varphi_0)_m)\in C([-h;0],\mathbb{R}^m)$, $k\leqslant m\leqslant n$,
or, more briefly, uniformly
$y$-bounded with $\varphi^{(m)}_0$-control, if for every number
$\alpha\geqslant 0$ there exists a number $\beta=\beta(\alpha)>\alpha$ such that for
every solution $x(t,t_0,\varphi_0)$, $\|\varphi^{(m)}_0\|_C\leqslant \alpha$, of equation $(2.1)$ the condition
$\|y_t(t_0,\varphi_0)\|_C\leqslant\beta$ is satisfied for all $t\geqslant t_0$, where $y_t(t_0,\varphi_0)=x_t(t_0,\varphi_0)^{(k)}$.

It is clear that if in Definition 3.2 we set $m=n$, then
we obtain the definition of uniform $y$-boundedness of solutions of equation $(2.1)$ given above. It is also worth noting that if in Definition 3.2 we set $k=n$, then we obtain the definition of uniform
$y$-boundedness with $\varphi^{(k)}_0$-control of solutions of equation $(2.1)$ introduced in the work \cite{Lap4}.

It is easy to see that the concept of uniform $y$-boundedness with $\varphi^{(m)}_0$-control of solutions of equation $(2.1)$ is equivalent to the fact that for every number $\alpha\geqslant 0$ there exists a $\beta=\beta(\alpha)>\alpha$ such that for every solution $x(t,t_0,\varphi_0)$, $\|\varphi^{(m)}_0\|_C\leqslant \alpha$, of equation $(2.1)$ the condition $\|y(t,t_0,\varphi_0)\|\leqslant\beta$ is satisfied for all $t\geqslant t_0$.

In what follows, for every element
$\xi=(\xi_1,\dots,\xi_l)\in\mathbb{R}^l$, $l\geqslant 1$, and any integers
$1\leqslant p\leqslant q\leqslant l$ the
notations $\mu=(\xi_1,\dots,\xi_p)\in\mathbb{R}^p$ and
$\gamma=(\xi_1,\dots,\xi_q)\in\mathbb{R}^q$ will be used.

Recall \cite{Lap1} (see also \cite{Lap2}) that solutions of the system of ordinary differential equations $\dot{\xi}=g(t,\xi)$,
where $\xi\in\mathbb{R}^l$, $l\geqslant 1$, and $g:\mathbb{R}^+\times\mathbb{R}^l\to\mathbb{R}^l$ is any continuous mapping, are said to be
uniformly bounded with respect to a part of the variables
$\mu=(\xi_1,\dots,\xi_p)$, $1\leqslant p\leqslant l$, with controlled part of initial conditions
$\gamma_0=((\xi_0)_1,\dots,(\xi_0)_q)$, $p\leqslant q\leqslant l$, or, more briefly, uniformly
$\mu$-bounded with $\gamma_0$-control, if for every number
$\nu\geqslant 0$ there exists a number $\lambda=\lambda(\nu)>\nu$ such that for
every solution $\xi(t,t_0,\xi_0)$, $\|\gamma_0\|\leqslant \nu$, of the system $\dot{\xi}=g(t,\xi)$ the condition
$\|\mu(t,t_0,\xi_0)\|\leqslant\lambda$ is satisfied for all $t\geqslant t_0$. Without loss of generality, we may assume that the function $\lambda(r)$ satisfies the condition $\lambda(r)\in\mathcal{K}$.
In the case when $p\leqslant q=l$, solutions of the above-mentioned system $\dot{\xi}=g(t,\xi)$ are said to be uniformly $\mu$-bounded \cite{Rum}. In the case when $p=q=l$, solutions of this system $\dot{\xi}=g(t,\xi)$ are said to be uniformly bounded \cite{Io2}. Methods for studying solutions of systems of ordinary differential equations for the above-mentioned types of uniform boundedness were developed, respectively, in the works \cite{Lap2}, \cite{Rum} and \cite{Io2}.

In what follows, by $\mathcal{K}$ we shall denote the set of all non-decreasing functions of the form $u:\mathbb{R}^+\to\mathbb{R}^+$.

Let us now formulate and prove a sufficient condition for uniform
$y$-boundedness with $\varphi_0^{(m)}$-control of solutions of equation $(2.1)$
in terms of vector Lyapunov-Krasovsky functionals.

{\bf Theorem 3.1}. Suppose that for equation $(2.1)$ there exist a vector Lyapunov-Krasovsky functional $V(t,\varphi)=(V_1(t,\varphi),\dots,V_l(t,\varphi))$,
$l\geqslant 1$,
with comparison system $(2.4)$ for equation $(2.1)$, functions $b(r),a(r)\in{\cal K}$, where $b(r)\to+\infty$ as $r\to+\infty$, integers $1\leqslant k\leqslant m\leqslant n$ and integers $1\leqslant
p\leqslant q\leqslant l$ for which the following conditions are satisfied:

1) $V_1(t,\varphi)\geqslant 0,\dots,V_q(t,\varphi)\geqslant 0$ for all $(t,\varphi)\in\mathbb{R}^+\times C([-h;0],\mathbb{R}^n)$;

2) $b(\|\varphi^{(k)}(0)\|)\leqslant\sum\limits_{i=1}^pV_i(t,\varphi)$ for all $(t,\varphi)\in\mathbb{R}^+\times C([-h;0],\mathbb{R}^n)$;

3) $\sum\limits_{i=1}^qV_i(t,\varphi)\leqslant a(\|\varphi^{(m)}\|_C)$ for all $(t,\varphi)\in\mathbb{R}^+\times C([-h;0],\mathbb{R}^n)$.

\noindent Moreover, suppose that the solutions of the comparison system $(2.4)$ for equation $(2.1)$ are uniformly $\mu$-bounded with
$\gamma_0$-control. Then the solutions of equation $(2.1)$ are uniformly
$y$-bounded with $\varphi_0^{(m)}$-control.

{\bf Proof}. It is required to show that for every number
$\alpha\geqslant 0$ and every solution $x(t,t_0,\varphi_0)$, where
$\|\varphi_0^{(m)}\|_C\leqslant\alpha$, of equation $(2.1)$ there exists a number
$\beta=\beta(\alpha)>\alpha$ such that the condition $\|y(t,t_0,\varphi_0)\|\leqslant\beta$
is satisfied for all $t\geqslant t_0$. Since $y(t,t_0,\varphi_0)=y_t(t_0,\varphi_0)(0)=x_t(t_0,\varphi_0)^{(k)}(0)$, using condition 2) of the theorem and
inequality $(2.5)$ we obtain for the solution $x(t,t_0,\varphi_0)$, where $\|\varphi_0^{(m)}\|_C\leqslant\alpha$, of equation $(2.1)$ and the upper
solution $\overline{\xi}(t,t_0,V(t_0,\varphi_0))$ of the comparison system $(2.4)$ the inequalities
$$b(\|y(t,t_0,\varphi_0)\|)\leqslant\sum_{i=1}^p
V_i(t,x_t(t_0,x_0))\leqslant
\sum_{i=1}^p\overline{\xi}_i(t,t_0,V(t_0,\varphi_0))$$ valid for all $t\geqslant t_0$.
Moreover, for all $t\geqslant t_0$
we have the obvious inequalities
$$\sum_{i=1}^p\overline{\xi}_i(t,t_0,V(t_0,\varphi_0))\leqslant
\sum_{i=1}^p|\,\overline{\xi}_i(t,t_0,V(t_0,\varphi_0))|\leqslant
p\cdot\|\,\overline{\mu}(t,t_0,V(t_0,\varphi_0))\|,$$ where
$\overline{\mu}(t,t_0,V(t_0,\varphi_0))=
(\overline{\xi}_1(t,t_0,V(t_0,\varphi_0)),\dots,\overline{\xi}_p(t,t_0,V(t_0,\varphi_0)))$.
Since by the condition of the theorem the solutions of the comparison system $(2.4)$ are
uniformly $\mu$-bounded with $\gamma_0$-control, for
the upper solution $\overline{\xi}(t,t_0,V(t_0,\varphi_0))$ of this system and the number
$\nu=\|V^\gamma(t_0,\varphi_0)\|$, where
$V^\gamma(t_0,\varphi_0)=(V_1(t_0,\varphi_0),\dots,V_q(t_0,\varphi_0))$, there is found
a number $\lambda(\nu)>\nu$ such that for all $t\geqslant t_0$
the condition
$\|\,\overline{\mu}(t,t_0,V(t_0,\varphi_0))\|\leqslant\lambda(\nu)$ is satisfied. From conditions
1) and 3) of the theorem we obtain
$$\|V^\gamma(t_0,\varphi_0)\|\leqslant\sum_{i=1}^q|\,V_i(t_0,\varphi_0)|=\sum_{i=1}^qV_i(t_0,\varphi_0)
\leqslant a(\|\varphi_0^{(m)}\|_C)\leqslant a(\alpha).$$
Since the function
$\lambda(r)$ is a non-decreasing function, using the inequalities
$\|V^\gamma(t_0,\varphi_0)\|\leqslant a(\alpha)$ and
$\|\,\overline{\mu}(t,t_0,V(t_0,\varphi_0))\|\leqslant\lambda(\|V^\gamma(t_0,\varphi_0)\|)$
we have
$$p\cdot\|\overline{\mu}(t,t_0,V(t_0,\varphi_0))\|\leqslant
p\cdot\lambda(\|V^\gamma(t_0,\varphi_0)\|)\leqslant
p\cdot\lambda(a(\alpha)).$$
From the inequalities obtained above
it follows that for all $t\geqslant t_0$ the inequality
$b(\|y(t,t_0,\varphi_0)\|)\leqslant p\cdot\lambda(a(\alpha))$ is valid. Using
now the condition $b(r)\to+\infty$ as $r\to+\infty$ and the fact that the number
$p$ is fixed, we choose a number $\beta>\alpha$, which depends on
$\alpha$ but does not depend on $t_0$, such that
$(p\cdot\lambda(a(\alpha)))\leqslant b(\beta)$. Using this inequality,
we obtain for all $t\geqslant t_0$ the inequality
$b(\|y(t,t_0,\varphi_0)\|)\leqslant b(\beta)$. Since the function $b(r)$ is
non-decreasing, from the last inequality we obtain
$||y(t,t_0,\varphi_0,)||\leqslant\beta$. Thus, it is shown that the solutions
of equation $(2.1)$ are uniformly $y$-bounded with
$\varphi_0^{(m)}$-control.~~~$\blacksquare$

If in the conditions of Theorem 3.1 we set $q=l$ or $p=q=l$, then we obtain, respectively, the following sufficient conditions
for uniform
$y$-boundedness with $\varphi_0^{(m)}$-control of solutions of equation $(2.1)$
in terms of vector Lyapunov-Krasovsky functionals.

{\bf Corollary 3.1}. Suppose that for equation $(2.1)$ there exist a vector Lyapunov-Krasovsky functional $V(t,\varphi)=(V_1(t,\varphi),\dots,V_l(t,\varphi))$,
$l\geqslant 1$,
with comparison system $(2.4)$ for equation $(2.1)$, functions $b(r),a(r)\in{\cal K}$, where $b(r)\to+\infty$ as $r\to+\infty$, integers $1\leqslant k\leqslant m\leqslant n$ and an integer $1\leqslant p\leqslant l$ for which the following conditions are satisfied:

1) $V_1(t,\varphi)\geqslant 0,\dots,V_l(t,\varphi)\geqslant 0$ for all $(t,\varphi)\in\mathbb{R}^+\times C([-h;0],\mathbb{R}^n)$;

2) $b(\|\varphi^{(k)}(0)\|)\leqslant\sum\limits_{i=1}^pV_i(t,\varphi)$ for all $(t,\varphi)\in\mathbb{R}^+\times C([-h;0],\mathbb{R}^n)$;

3) $\sum\limits_{i=1}^lV_i(t,\varphi)\leqslant a(\|\varphi^{(m)}\|_C)$ for all $(t,\varphi)\in\mathbb{R}^+\times C([-h;0],\mathbb{R}^n)$.

\noindent Moreover, suppose that the solutions of the comparison system $(2.4)$ for equation $(2.1)$ are uniformly $\mu$-bounded. Then the solutions of equation $(2.1)$ are uniformly
$y$-bounded with $\varphi_0^{(m)}$-control.~~~$\blacksquare$

{\bf Corollary 3.2}. Suppose that for equation $(2.1)$ there exist a vector Lyapunov-Krasovsky functional $V(t,\varphi)=(V_1(t,\varphi),\dots,V_l(t,\varphi))$,
$l\geqslant 1$,
with comparison system $(2.4)$ for equation $(2.1)$, functions $b(r),a(r)\in{\cal K}$, where $b(r)\to+\infty$ as $r\to+\infty$, and integers $1\leqslant k\leqslant m\leqslant n$ for which the following conditions are satisfied:

1) $V_1(t,\varphi)\geqslant 0,\dots,V_l(t,\varphi)\geqslant 0$ for all $(t,\varphi)\in\mathbb{R}^+\times C([-h;0],\mathbb{R}^n)$;

2) $b(\|\varphi^{(k)}(0)\|)\leqslant\sum\limits_{i=1}^lV_i(t,\varphi)\leqslant a(\|\varphi^{(m)}\|_C)$ for all $(t,\varphi)\in\mathbb{R}^+\times C([-h;0],\mathbb{R}^n)$.

\noindent Moreover, suppose that the solutions of the comparison system $(2.4)$ for equation $(2.1)$ are uniformly bounded. Then the solutions of equation $(2.1)$ are uniformly
$y$-bounded with $\varphi_0^{(m)}$-control.~~~$\blacksquare$

If in the conditions of Theorem 3.1 we set $m=n$, then we obtain the following statement, called, in view of its theoretical significance, also a theorem, which gives sufficient conditions for uniform $y$-boundedness
of solutions of equation $(2.1)$ in terms of vector Lyapunov-Krasovsky functionals.

{\bf Theorem 3.2}. Suppose that for equation $(2.1)$ there exist a vector Lyapunov-Krasovsky functional $V(t,\varphi)=(V_1(t,\varphi),\dots,V_l(t,\varphi))$,
$l\geqslant 1$,
with comparison system $(2.4)$ for equation $(2.1)$, functions $b(r),a(r)\in{\cal K}$, where $b(r)\to+\infty$ as $r\to+\infty$, an integer $1\leqslant k\leqslant n$ and integers $1\leqslant
p\leqslant q\leqslant l$ for which the following conditions are satisfied:

1) $V_1(t,\varphi)\geqslant 0,\dots,V_q(t,\varphi)\geqslant 0$ for all $(t,\varphi)\in\mathbb{R}^+\times C([-h;0],\mathbb{R}^n)$;

2) $b(\|\varphi^{(k)}(0)\|)\leqslant\sum\limits_{i=1}^pV_i(t,\varphi)$ for all $(t,\varphi)\in\mathbb{R}^+\times C([-h;0],\mathbb{R}^n)$;

3) $\sum\limits_{i=1}^qV_i(t,\varphi)\leqslant a(\|\varphi\|_C)$ for all $(t,\varphi)\in\mathbb{R}^+\times C([-h;0],\mathbb{R}^n)$.

\noindent Moreover, suppose that the solutions of the comparison system $(2.4)$ for equation $(2.1)$ are uniformly $\mu$-bounded with
$\gamma_0$-control. Then the solutions of equation $(2.1)$ are uniformly
$y$-bounded.~~~$\blacksquare$

If in the conditions of Theorem 3.2 we set $q=l$ or $p=q=l$, then we obtain, respectively, the following sufficient conditions
for uniform $y$-boundedness of solutions of equation $(2.1)$ in terms of vector Lyapunov-Krasovsky functionals.

{\bf Corollary 3.3}. Suppose that for equation $(2.1)$ there exist a vector Lyapunov-Krasovsky functional $V(t,\varphi)=(V_1(t,\varphi),\dots,V_l(t,\varphi))$,
$l\geqslant 1$,
with comparison system $(2.4)$ for equation $(2.1)$, functions $b(r),a(r)\in{\cal K}$, where $b(r)\to+\infty$ as $r\to+\infty$, an integer $1\leqslant k\leqslant n$ and an integer $1\leqslant
p\leqslant l$ for which the following conditions are satisfied:

1) $V_1(t,\varphi)\geqslant 0,\dots,V_l(t,\varphi)\geqslant 0$ for all $(t,\varphi)\in\mathbb{R}^+\times C([-h;0],\mathbb{R}^n)$;

2) $b(\|\varphi^{(k)}(0)\|)\leqslant\sum\limits_{i=1}^pV_i(t,\varphi)$ for all $(t,\varphi)\in\mathbb{R}^+\times C([-h;0],\mathbb{R}^n)$;

3) $\sum\limits_{i=1}^lV_i(t,\varphi)\leqslant a(\|\varphi\|_C)$ for all $(t,\varphi)\in\mathbb{R}^+\times C([-h;0],\mathbb{R}^n)$.

\noindent Moreover, suppose that the solutions of the comparison system $(2.4)$ for equation $(2.1)$ are uniformly $\mu$-bounded. Then the solutions of equation $(2.1)$ are uniformly
$y$-bounded.~$\blacksquare$

{\bf Corollary 3.4}. Suppose that for equation $(2.1)$ there exist a vector Lyapunov-Krasovsky functional $V(t,\varphi)=(V_1(t,\varphi),\dots,V_l(t,\varphi))$,
$l\geqslant 1$,
with comparison system $(2.4)$ for equation $(2.1)$, functions $b(r),a(r)\in{\cal K}$, where $b(r)\to+\infty$ as $r\to+\infty$, and an integer $1\leqslant k\leqslant n$ for which the following conditions are satisfied:

1) $V_1(t,\varphi)\geqslant 0,\dots,V_l(t,\varphi)\geqslant 0$ for all $(t,\varphi)\in\mathbb{R}^+\times C([-h;0],\mathbb{R}^n)$;

2) $b(\|\varphi^{(k)}(0)\|)\leqslant\sum\limits_{i=1}^lV_i(t,\varphi)\leqslant a(\|\varphi\|_C)$ for all $(t,\varphi)\in\mathbb{R}^+\times C([-h;0],\mathbb{R}^n)$.

\noindent Moreover, suppose that the solutions of the comparison system $(2.4)$ for equation $(2.1)$ are uniformly bounded. Then the solutions of equation $(2.1)$ are uniformly $y$-bounded.~~~$\blacksquare$

If in the conditions of Theorem 3.1 we set $k=m=n$, then we obtain the following statement, called, in view of its theoretical significance, also a theorem, which gives sufficient conditions for uniform boundedness
of solutions of equation $(2.1)$ in terms of vector Lyapunov-Krasovsky functionals.

{\bf Theorem 3.3}. Suppose that for equation $(2.1)$ there exist a vector Lyapunov-Krasovsky functional $V(t,\varphi)=(V_1(t,\varphi),\dots,V_l(t,\varphi))$,
$l\geqslant 1$,
with comparison system $(2.4)$ for equation $(2.1)$, functions $b(r),a(r)\in{\cal K}$, where $b(r)\to+\infty$ as $r\to+\infty$, and integers $1\leqslant
p\leqslant q\leqslant l$ for which the following conditions are satisfied:

1) $V_1(t,\varphi)\geqslant 0,\dots,V_q(t,\varphi)\geqslant 0$ for all $(t,\varphi)\in\mathbb{R}^+\times C([-h;0],\mathbb{R}^n)$;

2) $b(\|\varphi(0)\|)\leqslant\sum\limits_{i=1}^pV_i(t,\varphi)$ for all $(t,\varphi)\in\mathbb{R}^+\times C([-h;0],\mathbb{R}^n)$;

3) $\sum\limits_{i=1}^qV_i(t,\varphi)\leqslant a(\|\varphi\|_C)$ for all $(t,\varphi)\in\mathbb{R}^+\times C([-h;0],\mathbb{R}^n)$.

\noindent Moreover, suppose that the solutions of the comparison system $(2.4)$ for equation $(2.1)$ are uniformly $\mu$-bounded with
$\gamma_0$-control. Then the solutions of equation $(2.1)$ are uniformly
bounded.~~~$\blacksquare$

If in the conditions of Theorem 3.3 we set $q=l$ or $p=q=l$, then we obtain, respectively, the following sufficient conditions
for uniform boundedness of solutions of equation $(2.1)$ in terms of vector Lyapunov-Krasovsky functionals.

{\bf Corollary 3.5}. Suppose that for equation $(2.1)$ there exist a vector Lyapunov-Krasovsky functional $V(t,\varphi)=(V_1(t,\varphi),\dots,V_l(t,\varphi))$,
$l\geqslant 1$,
with comparison system $(2.4)$ for equation $(2.1)$, functions $b(r),a(r)\in{\cal K}$, where $b(r)\to+\infty$ as $r\to+\infty$, and an integer $1\leqslant
p\leqslant l$ for which the following conditions are satisfied:

1) $V_1(t,\varphi)\geqslant 0,\dots,V_l(t,\varphi)\geqslant 0$ for all $(t,\varphi)\in\mathbb{R}^+\times C([-h;0],\mathbb{R}^n)$;

2) $b(\|\varphi(0)\|)\leqslant\sum\limits_{i=1}^pV_i(t,\varphi)$ for all $(t,\varphi)\in\mathbb{R}^+\times C([-h;0],\mathbb{R}^n)$;

3) $\sum\limits_{i=1}^lV_i(t,\varphi)\leqslant a(\|\varphi\|_C)$ for all $(t,\varphi)\in\mathbb{R}^+\times C([-h;0],\mathbb{R}^n)$.

\noindent Moreover, suppose that the solutions of the comparison system $(2.4)$ for equation $(2.1)$ are uniformly $\mu$-bounded. Then the solutions of equation $(2.1)$ are uniformly
bounded.~~~$\blacksquare$

{\bf Corollary 3.6}. Suppose that for equation $(2.1)$ there exist a vector Lyapunov-Krasovsky functional $V(t,\varphi)=(V_1(t,\varphi),\dots,V_l(t,\varphi))$,
$l\geqslant 1$,
with comparison system $(2.4)$ for equation $(2.1)$ and functions $b(r),a(r)\in{\cal K}$, where $b(r)\to+\infty$ as $r\to+\infty$, for which the following conditions are satisfied:

1) $V_1(t,\varphi)\geqslant 0,\dots,V_l(t,\varphi)\geqslant 0$ for all $(t,\varphi)\in\mathbb{R}^+\times C([-h;0],\mathbb{R}^n)$;

2) $b(\|\varphi(0)\|)\leqslant\sum\limits_{i=1}^lV_i(t,\varphi)\leqslant a(\|\varphi\|_C)$ for all $(t,\varphi)\in\mathbb{R}^+\times C([-h;0],\mathbb{R}^n)$.

\noindent Moreover, suppose that the solutions of the comparison system $(2.4)$ for equation $(2.1)$ are uniformly bounded. Then the solutions of equation $(2.1)$ are uniformly bounded.~~~$\blacksquare$

It is important to note that Theorems 3.1\,--\,3.3 and Corollaries 3.1\,--\,3.6 make it possible, in the study of various types of uniform boundedness of solutions of the functional-differential equation $(2.1)$, to use methods for studying various types of uniform boundedness of solutions of systems of ordinary differential equations, since the comparison system $(2.4)$ for equation $(2.1)$ is a system of ordinary differential equations.

Consider now Theorem 3.1 in the case $l=1$, i.e., in the case when Lyapunov-Krasovsky functionals are taken as vector Lyapunov-Krasovsky functionals. Since in this case the comparison system $(2.4)$ for equation $(2.1)$ becomes the equation $\dot{\xi}=f(t,\xi)$, $\xi\in\mathbb{R}$, of comparison for equation $(2.1)$, we obtain the following statement, called, in view of its theoretical significance, also a theorem, which gives sufficient conditions
for uniform $y$-boundedness with $\varphi_0^{(m)}$-control of solutions of equation $(2.1)$ in terms of Lyapunov-Krasovsky functionals.

{\bf Theorem 3.4}. Suppose that for equation $(2.1)$ there exist a continuous functional $V:\mathbb{R}^+\times C([-h;0],\mathbb{R}^n)\to\mathbb{R}$ for which the derivative $\dot{V}_{(2.1)}(t,\varphi)$ along equation $(2.1)$ is also a continuous functional, a continuous
function $f(t,\xi)$ defined in $\mathbb{R}^+\times\mathbb{R}$, functions $b(r),a(r)\in{\cal K}$, where $b(r)\to+\infty$ as $r\to+\infty$, and integers $1\leqslant k\leqslant m\leqslant n$ such that the following conditions are satisfied:

1) $b(\|\|\varphi^{(k)}(0)\|)\leqslant V(t,\varphi)\leqslant a(\|\varphi^{(m)}\|_C)$ for all $(t,\varphi)\in\mathbb{R}^+\times C([-h;0],\mathbb{R}^n)$;

2) $\dot{V}_{(2.1)}(t,\varphi)\leqslant f(t,V(t,\varphi))$ for all
$(t,\varphi)\in\mathbb{R}^+\times C([-h;0],\mathbb{R}^n)$;

3) The solutions of the equation $\dot{\xi}=f(t,\xi)$ are uniformly bounded.

\noindent Then the solutions of equation $(2.1)$ are uniformly $y$-bounded
with $\varphi_0^{(m)}$-control.~~~$\blacksquare$

If in Theorem 3.4 we set $m=n$, then we obtain the following statement, called, in view of its theoretical significance, also a theorem, which gives sufficient conditions for uniform $y$-boundedness of solutions of equation $(2.1)$ in terms of Lyapunov-Krasovsky functionals.

{\bf Theorem 3.5}. Suppose that for equation $(2.1)$ there exist a continuous functional $V:\mathbb{R}^+\times C([-h;0],\mathbb{R}^n)\to\mathbb{R}$ for which the derivative $\dot{V}_{(2.1)}(t,\varphi)$ along equation $(2.1)$ is also a continuous functional, a continuous
function $f(t,\xi)$ defined in $\mathbb{R}^+\times\mathbb{R}$, functions $b(r),a(r)\in{\cal K}$, where $b(r)\to+\infty$ as $r\to+\infty$, and an integer $1\leqslant k\leqslant n$ such that the following conditions are satisfied:

1) $b(\|\|\varphi^{(k)}(0)\|)\leqslant V(t,\varphi)\leqslant a(\|\varphi\|_C)$ for all $(t,\varphi)\in\mathbb{R}^+\times C([-h;0],\mathbb{R}^n)$;

2) $\dot{V}_{(2.1)}(t,\varphi)\leqslant f(t,V(t,\varphi))$ for all
$(t,\varphi)\in\mathbb{R}^+\times C([-h;0],\mathbb{R}^n)$;

3) The solutions of the equation $\dot{\xi}=f(t,\xi)$ are uniformly bounded.

\noindent Then the solutions of equation $(2.1)$ are uniformly $y$-bounded.~~~$\blacksquare$

Finally, if in Theorem 3.4 we set $k=m=n$, then we obtain the following statement, called, in view of its theoretical significance, also a theorem, which gives sufficient conditions for uniform boundedness of solutions of equation $(2.1)$ in terms of Lyapunov-Krasovsky functionals.

{\bf Theorem 3.6}. Suppose that for equation $(2.1)$ there exist a continuous functional $V:\mathbb{R}^+\times C([-h;0],\mathbb{R}^n)\to\mathbb{R}$ for which the derivative $\dot{V}_{(2.1)}(t,\varphi)$ along equation $(2.1)$ is also a continuous functional, a continuous
function $f(t,\xi)$ defined in $\mathbb{R}^+\times\mathbb{R}$ and functions $b(r),a(r)\in{\cal K}$, where $b(r)\to+\infty$ as $r\to+\infty$, such that the following conditions are satisfied:

1) $b(\|\|\varphi(0)\|)\leqslant V(t,\varphi)\leqslant a(\|\varphi\|_C)$ for all $(t,\varphi)\in\mathbb{R}^+\times C([-h;0],\mathbb{R}^n)$;

2) $\dot{V}_{(2.1)}(t,\varphi)\leqslant f(t,V(t,\varphi))$ for all
$(t,\varphi)\in\mathbb{R}^+\times C([-h;0],\mathbb{R}^n)$;

3) The solutions of the equation $\dot{\xi}=f(t,\xi)$ are uniformly bounded.

\noindent Then the solutions of equation $(2.1)$ are uniformly bounded.~~~$\blacksquare$

Let us now proceed to the consideration of sufficient conditions for various kinds of uniform boundedness of solutions of equation $(2.1)$ in terms of higher derivatives of Lyapunov-Krasovsky functionals.

Let us recall \cite{Lap5} (see also \cite{Lap3}) that solutions of the ordinary differential equation of $l$-th order $\psi^{(l)}=g(t,\psi,\dot{\psi},\dots,\psi^{(l-1)})$, $\psi\in\mathbb{R}$, $l\geqslant 1$, where $g:\mathbb{R}^+\times\mathbb{R}^l\to\mathbb{R}^l$ is any continuous function, are said to be uniformly $(\psi,\dot{\psi},\dots,\psi^{(p-1)})$-bounded with $(\psi_0,\dot{\psi}_0,\dots,\psi^{(q-1)}_0)$-control, where
$1\leqslant p\leqslant q\leqslant l$, if the solutions of the corresponding system of ordinary differential equations of first order $\dot{\xi}=f(t,\xi)$, defined by formulas $(2.8)$, are uniformly $\mu$-bounded with $\gamma_0$-control.

A direct consequence of Theorems 2.2 and 3.1 is the following statement, called, in view of its theoretical significance, also a theorem, which gives sufficient conditions for uniform $y$-boundedness with $\varphi_0^{(m)}$-control of solutions of equation $(2.1)$
in terms of higher derivatives of Lyapunov-Krasovsky functionals.

{\bf Theorem 3.7}. Suppose that for equation $(2.1)$ there exist a Lyapunov-Krasovsky functional $U(t,\varphi)$ with an $l$-th order derivative along equation $(2.1)$ and comparison equation $(2.7)$ for equation $(2.1)$, functions $b(r),a(r)\in{\cal K}$, where $b(r)\to+\infty$ as $r\to+\infty$, integers $1\leqslant k\leqslant m\leqslant n$ and integers $1\leqslant
p\leqslant q\leqslant l$ for which the following conditions are satisfied:

1) $U^{(i)}_{(2.1)}(t,\varphi)\geqslant 0$ for every $0\leqslant i\leqslant
q-1$ and all $(t,\varphi)\in\mathbb{R}^+\times C([-h;0],\mathbb{R}^n)$;

2) $b(\|\varphi^{(k)}(0)\|)\leqslant\sum\limits_{i=1}^{p}U^{(i-1)}_{(2.1)}(t,\varphi)$ for all
$(t,\varphi)\in\mathbb{R}^+\times C([-h;0],\mathbb{R}^n)$;

3) $\sum\limits_{i=1}^{q}U^{(i-1)}_{(2.1)}(t,\varphi)\leqslant a(\|\varphi^{(m)}\|_C)$ for all $(t,\varphi)\in\mathbb{R}^+\times C([-h;0],\mathbb{R}^n)$.

\noindent Moreover, suppose that the solutions of the comparison equation $(2.7)$ for equation $(2.1)$
are uniformly $(\psi,\dot{\psi},\dots,\psi^{(p-1)})$-bounded with $(\psi_0,\dot{\psi}_0,\dots,\psi^{(q-1)}_0)$-control. Then the solutions of equation $(2.1)$ are uniformly
$y$-bounded with $\varphi_0^{(m)}$-control.~~~$\blacksquare$

If in the conditions of Theorem 3.7 we set $m=n$, then we obtain the following statement, called, in view of its theoretical significance, also a theorem, which gives sufficient conditions for uniform $y$-boundedness
of solutions of equation $(2.1)$ in terms of higher derivatives of Lyapunov-Krasovsky functionals.

{\bf Theorem 3.8}. Suppose that for equation $(2.1)$ there exist a Lyapunov-Krasovsky functional $U(t,\varphi)$ with an $l$-th order derivative along equation $(2.1)$ and comparison equation $(2.7)$ for equation $(2.1)$, functions $b(r),a(r)\in{\cal K}$, where $b(r)\to+\infty$ as $r\to+\infty$, an integer $1\leqslant k\leqslant n$ and integers $1\leqslant
p\leqslant q\leqslant l$ for which the following conditions are satisfied:

1) $U^{(i)}_{(2.1)}(t,\varphi)\geqslant 0$ for every $0\leqslant i\leqslant
q-1$ and all $(t,\varphi)\in\mathbb{R}^+\times C([-h;0],\mathbb{R}^n)$;

2) $b(\|\varphi^{(k)}(0)\|)\leqslant\sum\limits_{i=1}^{p}U^{(i-1)}_{(2.1)}(t,\varphi)$ for all
$(t,\varphi)\in\mathbb{R}^+\times C([-h;0],\mathbb{R}^n)$;

3) $\sum\limits_{i=1}^{q}U^{(i-1)}_{(2.1)}(t,\varphi)\leqslant a(\|\varphi\|_C)$ for all $(t,\varphi)\in\mathbb{R}^+\times C([-h;0],\mathbb{R}^n)$.

\noindent Moreover, suppose that the solutions of the comparison equation $(2.7)$ for equation $(2.1)$
are uniformly $(\psi,\dot{\psi},\dots,\psi^{(p-1)})$-bounded with $(\psi_0,\dot{\psi}_0,\dots,\psi^{(q-1)}_0)$-control. Then the solutions of equation $(2.1)$ are uniformly
$y$-bounded.~~~$\blacksquare$

If in the conditions of Theorem 3.7 we set $k=m=n$, then we obtain the following statement, called, in view of its theoretical significance, also a theorem, which gives sufficient conditions for uniform boundedness
of solutions of equation $(2.1)$ in terms of higher derivatives of Lyapunov-Krasovsky functionals.

{\bf Theorem 3.9}. Suppose that for equation $(2.1)$ there exist a Lyapunov-Krasovsky functional $U(t,\varphi)$ with an $l$-th order derivative along equation $(2.1)$ and comparison equation $(2.7)$ for equation $(2.1)$, functions $b(r),a(r)\in{\cal K}$, where $b(r)\to+\infty$ as $r\to+\infty$, and integers $1\leqslant
p\leqslant q\leqslant l$ for which the following conditions are satisfied:

1) $U^{(i)}_{(2.1)}(t,\varphi)\geqslant 0$ for every $0\leqslant i\leqslant
q-1$ and all $(t,\varphi)\in\mathbb{R}^+\times C([-h;0],\mathbb{R}^n)$;

2) $b(\|\varphi(0)\|)\leqslant\sum\limits_{i=1}^{p}U^{(i-1)}_{(2.1)}(t,\varphi)$ for all
$(t,\varphi)\in\mathbb{R}^+\times C([-h;0],\mathbb{R}^n)$;

3) $\sum\limits_{i=1}^{q}U^{(i-1)}_{(2.1)}(t,\varphi)\leqslant a(\|\varphi\|_C)$ for all $(t,\varphi)\in\mathbb{R}^+\times C([-h;0],\mathbb{R}^n)$.

\noindent Moreover, suppose that the solutions of the comparison equation $(2.7)$ for equation $(2.1)$
are uniformly $(\psi,\dot{\psi},\dots,\psi^{(p-1)})$-bounded with $(\psi_0,\dot{\psi}_0,\dots,\psi^{(q-1)}_0)$-control. Then the solutions of equation $(2.1)$ are uniformly
bounded.~~~$\blacksquare$

In conclusion of the section, let us give an example of the application of Corollary 3.6. Let a system of differential equations be given
$$\left\{\begin{array}{l}
\dot{x}_1(t)=-3x_1(t)-x_1^3(t)+x_1(t)\sin(x_2(t-h))\\
\dot{x}_2(t)=-3x_2(t)+x_2(t)\cos(x_1(t-h))\\
\end{array}\right.,\quad t\in\mathbb{R}^+,\eqno(3.1)$$
with constant delay $h>0$. Let us show, using Corollary 3.6, that the solutions of this system are uniformly bounded. Let us define for system $(3.1)$, for $l=2$, the vector Lyapunov-Krasovsky functional $V=(V_1(t,\varphi),V_2(t,\varphi))$, $\varphi=(\varphi_1,\varphi_2)\in C([-h;0],\mathbb{R}^2)$, by setting
$$V_1(t,\varphi)=\frac{1}{2}(\varphi_1^2(0)+\varphi_2^2(0)),$$
$$V_2(t,\varphi)=\frac{1}{2}(\varphi_1^2(0)+\varphi_2^2(0))+
\int_{-h}^0(\varphi_1^2(\vartheta)+\varphi_2^2(\vartheta)+\varphi_1^4(\vartheta))d\vartheta.$$
It is clear that condition 1) of Corollary 3.6 is satisfied. It is easy to see that for all $(t,\varphi)\in\mathbb{R}^+\times C([-h;0],\mathbb{R}^2)$ the double inequality
$$\|\varphi(0)\|^2\leqslant V_1(t,\varphi)+V_2(t,\varphi)\leqslant (1+h)\|\varphi\|_C^2+h\|\varphi\|_C^4$$
is valid, since $\varphi_1^2(0)+\varphi_2^2(0)=\|\varphi(0)\|^2\leqslant\|\varphi\|_C^2$, $\varphi_1^2(\vartheta)+\varphi_2^2(\vartheta)=\|\varphi(\vartheta)\|^2\leqslant\|\varphi\|_C^2$ and $\varphi_1^4(\vartheta)=\|\varphi_1(\vartheta)\|^4\leqslant\|\varphi\|_C^4$. Therefore, if as the
functions $b(r),a(r)\in{\cal K}$, where $b(r)\to+\infty$ as $r\to+\infty$, we take, respectively, the functions $b(r)=r^2$ and $a(r)=(1+h)r^2+hr^4$, then condition 2) of Corollary 3.6 is satisfied.
Let us now estimate
the derivative $(\dot{V}_1)_{(3.1)}(t,\varphi)=v_1'^+(\tau)|_{\tau=t}=v_1'^+(t)$ of the functional $V_1(t,\varphi)$ along system $(3.1)$, where
$v_1(\tau)=V_1(\tau,x_\tau)=\frac{1}{2}(x_1^2(\tau)+x_2^2(\tau))$
and $x(\tau,t,\varphi)=(x_1(\tau),x_2(\tau))$ is the solution of system $(3.1)$ with initial conditions $(t,\varphi)$. Since the function $v_1(\tau)$ is differentiable, for the derivative $(\dot{V}_1)_{(3.1)}(t,\varphi)=v_1'^+(t)=v_1'(t)$ we have
$$(\dot{V}_1)_{(3.1)}(t,\varphi)=v_1'(t)=x_1(t)\dot{x}_1(t)+x_2(t)\dot{x}_2(t)=$$
$$=x_1(t)(-3x_1(t)-x_1^3(t)+x_1(t)\sin(x_2(t-h)))+x_2(t)(-3x_2(t)+x_2(t)\cos(x_1(t-h)))=$$
$$=-3(x_1^2(t)+x_2^2(t))-x_1^4(t)+x_1^2(t)\sin(x_2(t-h))+x_2^2(t)\cos(x_1(t-h))\leqslant$$
$$\leqslant -2(x_1^2(t)+x_2^2(t))\leqslant -4V_1(t,\varphi)+2V_2(t,\varphi).$$
Let us now estimate
the derivative $(\dot{V}_2)_{(3.1)}(t,\varphi)=v_2'^+(\tau)|_{\tau=t}=v_2'^+(t)$ of the functional $V_2(t,\varphi)$ along system $(3.1)$, where
$$v_2(\tau)=V_2(\tau,x_\tau)=\frac{1}{2}(x_1^2(\tau)+x_2^2(\tau))+\int_{\tau-h}^\tau(x_1^2(s)+x_2^2(s)+x_1^4(s))ds$$
and $x(\tau,t,\varphi)=(x_1(\tau),x_2(\tau))$ is the solution of system $(3.1)$ with initial conditions $(t,\varphi)$.
Since the function $v_2(\tau)$ is differentiable, for the derivative $(\dot{V}_2)_{(3.1)}(t,\varphi)=v_2'^+(t)=v_2'(t)$ we have
$$(\dot{V}_2)_{(3.1)}(t,\varphi)=v_2'(t)=x_1(t)\dot{x}_1(t)+x_2(t)\dot{x}_2(t)+
\frac{d}{dt}\left(\int_{t-h}^t(x_1^2(s)+x_2^2(s)+x_1^4(s))ds\right)=$$
$$=x_1(t)(-3x_1(t)-x_1^3(t)+x_1(t)\sin(x_2(t-h)))+x_2(t)(-3x_2(t)+x_2(t)\cos(x_1(t-h)))\,+$$
$$+\,x_1^2(t)+x_2^2(t)+x_1^4(t)-x_1^2(t-h)-x_2^2(t-h)-x_1^4(t-h)\leqslant -(x_1^2(t)+x_2^2(t))=-2V_1(t,\varphi).$$
From the inequalities $(\dot{V}_1)_{(3.1)}(t,\varphi)\leqslant -4V_1(t,\varphi)+2V_2(t,\varphi)$ and $(\dot{V}_2)_{(3.1)}(t,\varphi)\leqslant -2V_1(t,\varphi)$ it follows that the comparison system for system $(3.1)$ is the system
$$\left\{\begin{array}{l}
\dot{\xi}_1=-4\xi_1+2\xi_2\\
\dot{\xi}_2=-2\xi_1\\
\end{array}\right.,\quad \xi=(\xi_1,\xi_2)\in\mathbb{R}^2.\eqno(3.2)$$
It is easy to verify that for the vector function $f(t,\xi)=(-4\xi_1+2\xi_2,-2\xi_1)$ the condition $f(t,\xi)\in W$ is satisfied. Let us show that the solutions of system $(3.2)$ are uniformly bounded. Consider for this system the Lyapunov function $W(t,\xi)=\xi_1^2+\xi_2^2$. It is clear that
$b(\|\xi\|)\leqslant  W(t,\xi)\leqslant a(\|\xi\|)$, where $b(r)=a(r)=r^2$, and $\dot{W}_{(3.2)}(t,\xi)=-8\xi_1^2\leqslant 0$. From this, by the Yoshizawa theorem on uniform boundedness of solutions of systems of ordinary differential equations \cite{Io2}, it follows that the solutions of system $(3.2)$ are uniformly bounded. Thus, for
system $(3.1)$ all conditions of Corollary 3.6 are satisfied and, consequently, the solutions of system $(3.1)$
are uniformly bounded.
\vspace{0.5cm}

\centerline{\bf 4. Vector Lyapunov-Krasovsky functionals, higher derivatives}
\centerline{\bf of Lyapunov-Krasovsky functionals and ultimate-uniform boundedness}
\centerline{\bf of solutions of functional-differential equations of retarded type}
\vspace{0.5cm}

Let us first recall \cite{Io1} that solutions of equation $(2.1)$ are said to be ultimate-uniformly bounded if for equation $(2.1)$ there is found a number $B>0$, and for every
$\alpha\geqslant 0$ there exists a number $T=T(\alpha)\geqslant 0$ such that for
every solution $x(t,t_0,\varphi_0)$, where $\|\varphi_0\|_C\leqslant \alpha$, of equation $(2.1)$ the condition
$\|x_t(t_0,\varphi_0)\|_C<B$ is satisfied for all $t\geqslant t_0+T$.

It is easy to see that the concept of ultimate-uniform boundedness of solutions of equation $(2.1)$ is equivalent to the fact that for equation $(2.1)$ there is found a number $B>0$, and for every
$\alpha\geqslant 0$ there exists a number $T=T(\alpha)\geqslant 0$ such that for
every solution $x(t,t_0,\varphi_0)$, where $\|\varphi_0\|_C\leqslant \alpha$, of equation $(2.1)$ the condition
$\|x(t,t_0,\varphi_0)\|<B$ is satisfied for all $t\geqslant t_0+T-h$.

Let us now recall \cite{Lap4} that solutions of equation $(2.1)$ are said to be ultimate-uniformly bounded with respect to a part of the variables
$y=(x_1,\dots,x_k)$, $1\leqslant k\leqslant n$, with partial control of prehistories $\varphi^{(k)}_0\in C([-h;0],\mathbb{R}^k)$ or, more briefly, ultimate-uniformly $y$-bounded with $\varphi^{(k)}_0$-control, if for equation $(2.1)$ there is found a number $B>0$, and for every
$\alpha\geqslant 0$ there exists a number $T=T(\alpha)\geqslant 0$ such that for
every solution $x(t,t_0,\varphi_0)$, where $\|\varphi^{(k)}_0\|_C\leqslant \alpha$, of equation $(2.1)$ the condition
$\|y_t(t_0,\varphi_0)\|_C<B$ is satisfied for all $t\geqslant t_0+T$.

Obviously, for $k=n$ the concept of ultimate-uniform $y$-boundedness with $\varphi^{(k)}_0$-control of solutions of equation $(2.1)$ becomes the concept of ultimate-uniform boundedness of solutions of this equation.

It is easy to see that the concept of ultimate-uniform $y$-boundedness with $\varphi^{(k)}_0$-control of solutions of equation $(2.1)$ is equivalent to the fact that for equation $(2.1)$ there is found a number $B>0$, and for every
$\alpha\geqslant 0$ there exists a number $T=T(\alpha)\geqslant 0$ such that for
every solution $x(t,t_0,\varphi_0)$, where $\|\varphi_0\|_C\leqslant \alpha$, of equation $(2.1)$ the condition
$\|y(t,t_0,\varphi_0)\|<B$ is satisfied for all $t\geqslant t_0+T-h$.

Recall \cite{Lap} (see also \cite{Lap3}) that solutions of the system of ordinary differential equations $\dot{\xi}=g(t,\xi)$,
where $\xi\in\mathbb{R}^l$, $l\geqslant 1$, and $g:\mathbb{R}^+\times\mathbb{R}^l\to\mathbb{R}^l$ is any continuous mapping, are said to be
ultimate-uniformly bounded with respect to a part of the variables
$\mu=(\xi_1,\dots,\xi_p)$ with controlled part of initial conditions
$\mu_0=((\xi_0)_1,\dots,(\xi_0)_p)$ or, more briefly, ultimate-uniformly
$\mu$-bounded with $\mu_0$-control, if for the system $\dot{\xi}=g(t,\xi)$ there is found a number $D>0$, and for every
$\nu\geqslant 0$ there exists a number $\tau=\tau(\nu)\geqslant 0$ such that for
every solution $\xi(t,t_0,\xi_0)$, where $\|\mu_0\|\leqslant \nu$, of the system $\dot{\xi}=g(t,\xi)$ the condition
$\|\mu(t,t_0,\xi_0)\|<D$ is satisfied for all $t\geqslant t_0+\tau$. In the case when $p=l$, solutions of the system $\dot{\xi}=g(t,\xi)$ are said to be ultimate-uniformly bounded \cite{Io2}. Methods for studying solutions of systems of ordinary differential equations for the above-mentioned types of ultimate-uniform boundedness were developed, respectively, in the works \cite{Lap} and \cite{Io2}.

Let us now formulate and prove a sufficient condition for ultimate-uniform
$y$-boundedness with $\varphi_0^{(k)}$-control of solutions of equation $(2.1)$
in terms of vector Lyapunov-Krasovsky functionals.

{\bf Theorem 4.1}. Suppose that for equation $(2.1)$ there exist a vector Lyapunov-Krasovsky functional $V(t,\varphi)=(V_1(t,\varphi),\dots,V_l(t,\varphi))$,
$l\geqslant 1$,
with comparison system $(2.4)$ for equation $(2.1)$, functions $b(r),a(r)\in{\cal K}$, where $b(r)\to+\infty$ as $r\to+\infty$, and integers $1\leqslant k\leqslant n$ and $1\leqslant p\leqslant l$ for which the following conditions are satisfied:

1) $V_1(t,\varphi)\geqslant 0,\dots,V_p(t,\varphi)\geqslant 0$ for all $(t,\varphi)\in\mathbb{R}^+\times C([-h;0],\mathbb{R}^n)$;

2) $b(\|\varphi^{(k)}(0)\|)\leqslant\sum\limits_{i=1}^pV_i(t,\varphi)\leqslant a(\|\varphi^{(k)}\|_C)$ for all $(t,\varphi)\in\mathbb{R}^+\times C([-h;0],\mathbb{R}^n)$.

\noindent Moreover, suppose that the solutions of the comparison system $(2.4)$ for equation $(2.1)$ are ultimate-uniformly $\mu$-bounded with
$\mu_0$-control. Then the solutions of equation $(2.1)$ are ultimate-uniformly
$y$-bounded with $\varphi_0^{(k)}$-control.

{\bf Proof}. It is required to show that for equation $(2.1)$ there is found a number $B>0$, and for every
$\alpha\geqslant 0$ there exists a number $T=T(\alpha)\geqslant 0$ such that for
every solution $x(t,t_0,\varphi_0)$, where $\|\varphi^{(k)}_0\|_C\leqslant \alpha$, of equation $(2.1)$ the condition
$\|y_t(t_0,\varphi_0)\|_C<B$ is satisfied for all $t\geqslant t_0+T$. Since $y(t,t_0,\varphi_0)=y_t(t_0,\varphi_0)(0)=x_t(t_0,\varphi_0)^{(k)}(0)$, using the left inequality from condition 2) of the theorem and
inequality $(2.5)$ we obtain for the solution $x(t,t_0,\varphi_0)$, where $\|\varphi_0^{(k)}\|_C\leqslant\alpha$, of equation $(2.1)$ and the upper
solution $\overline{\xi}(t,t_0,V(t_0,\varphi_0))$ of the comparison system $(2.4)$ the inequalities $$b(\|y(t,t_0,\varphi_0)\|)\leqslant\sum_{i=1}^p
V_i(t,x_t(t_0,x_0))\leqslant
\sum_{i=1}^p\overline{\xi}_i(t,t_0,V(t_0,\varphi_0))$$ valid for all $t\geqslant t_0$.
Moreover, for all $t\geqslant t_0$
we have the obvious inequalities
$$\sum_{i=1}^p\overline{\xi}_i(t,t_0,V(t_0,\varphi_0))\leqslant
\sum_{i=1}^p|\,\overline{\xi}_i(t,t_0,V(t_0,\varphi_0))|\leqslant
p\cdot\|\,\overline{\mu}(t,t_0,V(t_0,\varphi_0))\|,$$ where
$\overline{\mu}(t,t_0,V(t_0,\varphi_0))=
(\overline{\xi}_1(t,t_0,V(t_0,\varphi_0)),\dots,\overline{\xi}_p(t,t_0,V(t_0,\varphi_0)))$.
Since by the condition of the theorem the solutions of the comparison system $(2.4)$ are ultimate-uniformly $\mu$-bounded with $\mu_0$-control, for this system there exists a
number $D>0$ such that for the upper solution $\overline{\xi}(t,t_0,V(t_0,\varphi_0))$ of system $(2.4)$ and the number
$\nu=\|V^\mu(t_0,\varphi_0)\|$, where $V^\mu(t_0,\varphi_0)=(V_1(t_0,\varphi_0),\dots,V_p(t_0,\varphi_0))$, there is found
a number $\tau=\tau(\nu)\geqslant 0$ such that the condition $\|\overline{\mu}(t,t_0,V(t_0,\varphi_0))\|<D$ is satisfied for all $t\geqslant t_0+\tau$.
From this and the inequalities indicated above we have, for all $t\geqslant t_0+\tau$, the inequality
$b(\|y(t,t_0,\varphi_0)\|)<p\cdot D$. Using condition 1) and the right inequality from condition 2), we obtain
$$\|V^\mu(t_0,\varphi_0)\|\leqslant\sum_{i=1}^p|\,V_i(t_0,\varphi_0)|=\sum_{i=1}^pV_i(t_0,\varphi_0)
\leqslant a(\|\varphi_0^{(k)}\|_C)\leqslant a(\alpha).$$
Since the function $\tau(r)$ is a non-decreasing function, we have $\tau(\|V^\mu(t_0,\varphi_0)\|)\leqslant\tau(a(\alpha))$.
If we define the required number $T=T(\alpha)$ by the equality $T(\alpha)=\tau(a(\alpha))+h$, then for this number $T$ we have $b(\|y(t,t_0,\varphi_0)\|)<p\cdot D$
for all $t\geqslant t_0+T-h$. Using now the condition $b(r)\to+\infty$ as
$r\to+\infty$ and the fact that the numbers $p$ and $D$ are fixed, we choose a number $B>0$
such that $p\cdot D<b(B)$. From this, for all $t\geqslant t_0+T-h$, we obtain the inequality $b(y(t,t_0,\varphi_0))<b(B)$, where $B$ does not depend on the solution $x(t,t_0,\varphi_0)$ but depends only on equation $(2.1)$. Since
the function $b(r)$ is non-decreasing, from the last inequality
we obtain $\|y(t,t_0,\varphi_0)\|<B$ for all $t\geqslant t_0+T-h$, which is equivalent to the validity of the inequality $\|y_t(t_0,\varphi_0)\|_C<B$ for all $t\geqslant t_0+T$.
Thus, it is shown that the solutions of equation $(2.1)$ are ultimate-uniformly $y$-bounded with
$\varphi_0^{(k)}$-control.~~~$\blacksquare$

If in the conditions of Theorem 4.1 we set $p=l$, then we obtain the following sufficient condition for ultimate-uniform $y$-boundedness with $\varphi_0^{(k)}$-control of solutions of equation $(2.1)$
in terms of vector Lyapunov-Krasovsky functionals.

{\bf Corollary 4.1}. Suppose that for equation $(2.1)$ there exist a vector Lyapunov-Krasovsky functional $V(t,\varphi)=(V_1(t,\varphi),\dots,V_l(t,\varphi))$,
$l\geqslant 1$,
with comparison system $(2.4)$ for equation $(2.1)$, functions $b(r),a(r)\in{\cal K}$, where $b(r)\to+\infty$ as $r\to+\infty$, and an integer $1\leqslant k\leqslant n$ for which the following conditions are satisfied:

1) $V_1(t,\varphi)\geqslant 0,\dots,V_l(t,\varphi)\geqslant 0$ for all $(t,\varphi)\in\mathbb{R}^+\times C([-h;0],\mathbb{R}^n)$;

2) $b(\|\varphi^{(k)}(0)\|)\leqslant\sum\limits_{i=1}^lV_i(t,\varphi)\leqslant a(\|\varphi^{(k)}\|_C)$ for all $(t,\varphi)\in\mathbb{R}^+\times C([-h;0],\mathbb{R}^n)$.

\noindent Moreover, suppose that the solutions of the comparison system $(2.4)$ for equation $(2.1)$ are ultimate-uniformly bounded. Then the solutions of equation $(2.1)$ are ultimate-uniformly
$y$-bounded with $\varphi_0^{(k)}$-control.~~~$\blacksquare$

If in the conditions of Theorem 4.1 we set $k=n$, then we obtain the following statement, called, in view of its theoretical significance, also a theorem, which gives sufficient conditions for ultimate-uniform boundedness of solutions of equation $(2.1)$
in terms of vector Lyapunov-Krasovsky functionals.

{\bf Theorem 4.2}. Suppose that for equation $(2.1)$ there exist a vector Lyapunov-Krasovsky functional $V(t,\varphi)=(V_1(t,\varphi),\dots,V_l(t,\varphi))$,
$l\geqslant 1$,
with comparison system $(2.4)$ for equation $(2.1)$, functions $b(r),a(r)\in{\cal K}$, where $b(r)\to+\infty$ as $r\to+\infty$, and an integer $1\leqslant p\leqslant l$ for which the following conditions are satisfied:

1) $V_1(t,\varphi)\geqslant 0,\dots,V_p(t,\varphi)\geqslant 0$ for all $(t,\varphi)\in\mathbb{R}^+\times C([-h;0],\mathbb{R}^n)$;

2) $b(\|\varphi(0)\|)\leqslant\sum\limits_{i=1}^pV_i(t,\varphi)\leqslant a(\|\varphi\|_C)$ for all $(t,\varphi)\in\mathbb{R}^+\times C([-h;0],\mathbb{R}^n)$.

\noindent Moreover, suppose that the solutions of the comparison system $(2.4)$ for equation $(2.1)$ are ultimate-uniformly $\mu$-bounded with
$\mu_0$-control. Then the solutions of equation $(2.1)$ are ultimate-uniformly
bounded.~~~$\blacksquare$

If in the conditions of Theorem 4.2 we set $p=l$, then we obtain the following sufficient condition for ultimate-uniform boundedness of solutions of equation $(2.1)$ in terms of vector Lyapunov-Krasovsky functionals.

{\bf Corollary 4.2}. Suppose that for equation $(2.1)$ there exist a vector Lyapunov-Krasovsky functional $V(t,\varphi)=(V_1(t,\varphi),\dots,V_l(t,\varphi))$,
$l\geqslant 1$,
with comparison system $(2.4)$ for equation $(2.1)$ and functions $b(r),a(r)\in{\cal K}$, where $b(r)\to+\infty$ as $r\to+\infty$, for which the following conditions are satisfied:

1) $V_1(t,\varphi)\geqslant 0,\dots,V_l(t,\varphi)\geqslant 0$ for all $(t,\varphi)\in\mathbb{R}^+\times C([-h;0],\mathbb{R}^n)$;

2) $b(\|\varphi(0)\|)\leqslant\sum\limits_{i=1}^lV_i(t,\varphi)\leqslant a(\|\varphi\|_C)$ for all $(t,\varphi)\in\mathbb{R}^+\times C([-h;0],\mathbb{R}^n)$.

\noindent Moreover, suppose that the solutions of the comparison system $(2.4)$ for equation $(2.1)$ are ultimate-uniformly bounded. Then the solutions of equation $(2.1)$ are ultimate-uniformly bounded.~~~$\blacksquare$

It is important to note that Theorems 4.1, 4.2 and Corollaries 4.1, 4.2 make it possible, in the study of various types of ultimate-uniform boundedness of solutions of the functional-differential equation $(2.1)$, to use methods for studying various types of ultimate-uniform boundedness of solutions of systems of ordinary differential equations, since the comparison system $(2.4)$ for equation $(2.1)$ is a system of ordinary differential equations.

Consider now Theorem 4.1 in the case $l=1$, i.e., in the case when Lyapunov-Krasovsky functionals are taken as vector Lyapunov-Krasovsky functionals. Since in this case the comparison system $(3.2)$ for equation $(2.1)$ becomes the equation $\dot{\xi}=f(t,\xi)$, $\xi\in\mathbb{R}$, of comparison for equation $(2.1)$, we obtain the following statement, called, in view of its theoretical significance, also a theorem, which gives sufficient conditions
for ultimate-uniform $y$-boundedness with $\varphi_0^{(k)}$-control of solutions of equation $(2.1)$ in terms of Lyapunov-Krasovsky functionals.

{\bf Teorem 4.3}. Suppose that for equation $(2.1)$ there exist a continuous functional $V:\mathbb{R}^+\times C([-h;0],\mathbb{R}^n)\to\mathbb{R}$ for which the derivative $\dot{V}(t,\varphi)$ along equation $(2.1)$ is also a continuous functional, a continuous
function $f(t,\xi)$ defined in $\mathbb{R}^+\times\mathbb{R}$, functions $b(r),a(r)\in{\cal K}$, where $b(r)\to+\infty$ as $r\to+\infty$, and an integer $1\leqslant k\leqslant n$ such that the following conditions are satisfied:

1) $b(\|\|\varphi^{(k)}(0)\|)\leqslant V(t,\varphi)\leqslant a(\|\varphi^{(k)}\|_C)$ for all $(t,\varphi)\in\mathbb{R}^+\times C([-h;0],\mathbb{R}^n)$;

2) $\dot{V}(t,\varphi)\leqslant f(t,V(t,\varphi))$ for all
$(t,\varphi)\in\mathbb{R}^+\times C([-h;0],\mathbb{R}^n)$;

3) The solutions of the equation $\dot{\xi}=f(t,\xi)$ are ultimate-uniformly bounded.

\noindent Then the solutions of equation $(2.1)$ are ultimate-uniformly $y$-bounded
with $\varphi_0^{(k)}$-control.~~~$\blacksquare$

If in the conditions of Theorem 4.3 we set $k=n$, then we obtain the following sufficient condition for ultimate-uniform boundedness of solutions of equation $(2.1)$ in terms of Lyapunov-Krasovsky functionals.

{\bf Corollary 4.3}. Suppose that for equation $(2.1)$ there exist a continuous functional $V:\mathbb{R}^+\times C([-h;0],\mathbb{R}^n)\to\mathbb{R}$ for which the derivative $\dot{V}(t,\varphi)$ along equation $(2.1)$ is also a continuous functional, a continuous
function $f(t,\xi)$ defined in $\mathbb{R}^+\times\mathbb{R}$, and functions $b(r),a(r)\in{\cal K}$, where $b(r)\to+\infty$ as $r\to+\infty$, such that the following conditions are satisfied:

1) $b(\|\|\varphi(0)\|)\leqslant V(t,\varphi)\leqslant a(\|\varphi\|_C)$ for all $(t,\varphi)\in\mathbb{R}^+\times C([-h;0],\mathbb{R}^n)$;

2) $\dot{V}(t,\varphi)\leqslant f(t,V(t,\varphi))$ for all
$(t,\varphi)\in\mathbb{R}^+\times C([-h;0],\mathbb{R}^n)$;

3) The solutions of the equation $\dot{\xi}=f(t,\xi)$ are ultimate-uniformly bounded.

\noindent Then the solutions of equation $(2.1)$ are ultimate-uniformly bounded.~~~$\blacksquare$

Let us now proceed to the consideration of sufficient conditions for various kinds of ultimate-uniform boundedness of solutions of equation $(2.1)$ in terms of higher derivatives of Lyapunov-Krasovsky functionals.

Let us recall \cite{Lap6} (see also \cite{Lap3}) that solutions of the ordinary differential equation of $l$-th order $\psi^{(l)}=g(t,\psi,\dot{\psi},\dots,\psi^{(l-1)})$, $\psi\in\mathbb{R}$, $l\geqslant 1$, where $g:\mathbb{R}^+\times\mathbb{R}^l\to\mathbb{R}^l$ is any continuous function, are said to be ultimate-uniformly $(\psi,\dot{\psi},\dots,\psi^{(p-1)})$-bounded with $(\psi_0,\dot{\psi}_0,\dots,\psi^{(p-1)}_0)$-control, where
$1\leqslant p\leqslant l$, if the solutions of the corresponding system of ordinary differential equations of first order $\dot{\xi}=f(t,\xi)$, defined by formulas $(2.8)$, are ultimate-uniformly $\mu$-bounded with $\mu_0$-control.

A direct consequence of Theorems 2.2 and 4.1 is the following statement, called, in view of its theoretical significance, also a theorem, which gives sufficient conditions for ultimate-uniform $y$-boundedness with $\varphi_0^{(k)}$-control of solutions of equation $(2.1)$
in terms of higher derivatives of Lyapunov-Krasovsky functionals.

{\bf Theorem 4.4}. Suppose that for equation $(2.1)$ there exist a Lyapunov-Krasovsky functional $U(t,\varphi)$ with an $l$-th order derivative along equation $(2.1)$ and comparison equation $(2.7)$ for equation $(2.1)$, functions $b(r),a(r)\in{\cal K}$, where $b(r)\to+\infty$ as $r\to+\infty$, and integers $1\leqslant k\leqslant n$ and $1\leqslant
p\leqslant l$ for which the following conditions are satisfied:

\noindent 1) $U^{(i)}_{(2.1)}(t,\varphi)\geqslant 0$ for every $0\leqslant i\leqslant
p-1$ and all $(t,\varphi)\in\mathbb{R}^+\times C([-h;0],\mathbb{R}^n)$;

\noindent 2) $b(\|\varphi^{(k)}(0)\|)\leqslant\sum\limits_{i=1}^{p}U^{(i-1)}_{(2.1)}(t,\varphi)\leqslant a(\|\varphi^{(k)}\|_C)$ for all
$(t,\varphi)\in\mathbb{R}^+\times C([-h;0],\mathbb{R}^n)$.

\noindent Moreover, suppose that the solutions of the comparison equation $(2.7)$ for equation $(2.1)$
are ultimate-uniformly $(\psi,\dot{\psi},\dots,\psi^{(p-1)})$-bounded with $(\psi_0,\dot{\psi}_0,\dots,\psi^{(p-1)}_0)$-control. Then the solutions of equation $(2.1)$ are ultimate-uniformly
$y$-bounded with $\varphi_0^{(k)}$-control.~~~$\blacksquare$

If in the conditions of Theorem 4.4 we set $k=n$, then we obtain the following statement, called, in view of its theoretical significance, also a theorem, which gives sufficient conditions for ultimate-uniform boundedness
of solutions of equation $(2.1)$ in terms of higher derivatives of Lyapunov-Krasovsky functionals.

{\bf Theorem 4.5}. Suppose that for equation $(2.1)$ there exist a Lyapunov-Krasovsky functional $U(t,\varphi)$ with an $l$-th order derivative along equation $(2.1)$ and comparison equation $(2.7)$ for equation $(2.1)$, functions $b(r),a(r)\in{\cal K}$, where $b(r)\to+\infty$ as $r\to+\infty$, and integer $1\leqslant
p\leqslant l$ for which the following conditions are satisfied:

\noindent 1) $U^{(i)}_{(2.1)}(t,\varphi)\geqslant 0$ for every $0\leqslant i\leqslant
p-1$ and all $(t,\varphi)\in\mathbb{R}^+\times C([-h;0],\mathbb{R}^n)$;

\noindent 2) $b(\|\varphi(0)\|)\leqslant\sum\limits_{i=1}^{p}U^{(i-1)}_{(2.1)}(t,\varphi)\leqslant a(\|\varphi\|_C)$ for all
$(t,\varphi)\in\mathbb{R}^+\times C([-h;0],\mathbb{R}^n)$.

\noindent Moreover, suppose that the solutions of the comparison equation $(2.7)$ for equation $(2.1)$
are ultimate-uniformly $(\psi,\dot{\psi},\dots,\psi^{(p-1)})$-bounded with $(\psi_0,\dot{\psi}_0,\dots,\psi^{(p-1)}_0)$-control. Then the solutions of equation $(2.1)$ are ultimate-uniformly
bounded.~~~$\blacksquare$

In conclusion of the section, let us give an example of the application of Theorem 4.1. Let a system of differential equations be given
$$\left\{\begin{array}{l}
\dot{x}_1(t)=-3x_1(t)+x_1(t)\cos(x_2(t-h))\\
\dot{x}_2(t)=-3x_2(t)+x_2(t)\sin(x_1(t-h))\\
\end{array}\right.,\quad t\in\mathbb{R}^+,\eqno(4.1)$$
with constant delay $h>0$. Let us show, using Theorem 4.1, that the solutions of this system are ultimate-uniformly $y$-bounded with $\varphi_0^{(1)}$-control, where $y=(x_1)$ and $\varphi_0^{(1)}=((\varphi_0)_1)$. Let us define for system $(4.1)$, for $l=2$, the vector Lyapunov-Krasovsky functional $V=(V_1(t,\varphi),V_2(t,\varphi))$, $\varphi=(\varphi_1,\varphi_2)\in C([-h;0],\mathbb{R}^2)$, by setting
$$V_1(t,\varphi)=\frac{1}{2}\,\varphi_1^2(0),\quad V_2(t,\varphi)=\frac{1}{2}(\varphi_1^2(0)+\varphi_2^2(0))+
\int_{-h}^0(\varphi_1^2(\vartheta)+\varphi_2^2(\vartheta))d\vartheta.$$
It is clear that for $p=1$ condition 1) of Theorem 4.1 is satisfied. It is easy to see that for all $(t,\varphi)\in\mathbb{R}^+\times C([-h;0],\mathbb{R}^2)$ the double inequality
$$\frac{1}{2}\|\varphi^{(1)}(0)\|^2\leqslant V_1(t,\varphi)\leqslant\frac{1}{2}\|\varphi^{(1)}\|_C^2$$
is valid, since $\varphi^{(1)}(0)=\varphi_1(0)$, $\|\varphi^{(1)}(0)\|^2=(\varphi^{(1)}(0))^2=\varphi_1^2(0)$ and $\|\varphi^{(1)}(0)\|^2\leqslant\|\varphi^{(1)}\|_C^2$. Therefore, if as the
functions $b(r),a(r)\in{\cal K}$ we take, respectively, the functions $b(r)=a(r)=(1/2)r^2$, then for $p=1$ condition 2) of Theorem 4.1 is satisfied.
Let us now estimate
the derivative $(\dot{V}_1)_{(4.1)}(t,\varphi)=v_1'^+(\tau)|_{\tau=t}=v_1'^+(t)$ of the functional $V_1(t,\varphi)$ along system $(4.1)$, where
$v_1(\tau)=V_1(\tau,x_\tau(t,\varphi))=(1/2)x_1^2(\tau)$
and $x(\tau,t,\varphi)=(x_1(\tau),x_2(\tau))$ is the solution of system $(4.1)$ with initial conditions $(t,\varphi)$.
Since the function $v_1(\tau)$ is differentiable, for the derivative $(\dot{V}_1)_{(4.1)}(t,\varphi)=v_1'^+(t)=v_1'(t)$ we have
$$(\dot{V}_1)_{(4.1)}(t,\varphi)=v_1'(t)=x_1(t)\dot{x}_1(t)=x_1(t)(-3x_1(t)+x_1(t)\cos(x_2(t-h)))=$$
$$=-3x_1^2(t)+x_1^2(t)\cos(x_2(t-h))\leqslant -2x_1^2(t)=-2x_1^2(\tau)|_{\tau=t}=$$
$$=-4V_1(\tau,x_\tau(t,\varphi))|_{\tau=t}=-4V_1(t,\varphi).$$
Let us now estimate
the derivative $(\dot{V}_2)_{(4.1)}(t,\varphi)=v_2'^+(\tau)|_{\tau=t}=v_2'^+(t)$ of the functional $V_2(t,\varphi)$ along system $(4.1)$, where
$$v_2(\tau)=V_2(\tau,x_\tau(t,\varphi))=\frac{1}{2}(x_1^2(\tau)+x_2^2(\tau))+\int_{\tau-h}^\tau(x_1^2(s)+x_2^2(s))ds$$
and $x(\tau,t,\varphi)=(x_1(\tau),x_2(\tau))$ is the solution of system $(4.1)$ with initial conditions $(t,\varphi)$.
Since the function $v_2(\tau)$ is differentiable, for the derivative $(\dot{V}_2)_{(4.1)}(t,\varphi)=v_2'^+(t)=v_2'(t)$ we have
$$(\dot{V}_2)_{(4.1)}(t,\varphi)=v_2'(t)=x_1(t)\dot{x}_1(t)+x_2(t)\dot{x}_2(t)+
\frac{d}{dt}\left(\int_{t-h}^t(x_1^2(s)+x_2^2(s))ds\right)=$$
$$=x_1(t)(-3x_1(t)+x_1(t)\cos(x_2(t-h)))+x_2(t)(-3x_2(t)+x_2(t)\sin(x_1(t-h)))\,+$$
$$+\,x_1^2(t)+x_2^2(t)-x_1^2(t-h)-x_2^2(t-h)\leqslant -x_1^2(t)-x_2^2(t)\leqslant-x_1^2(t)=$$
$$=-x_1^2(\tau)|_{\tau=t}=-2V_1(\tau,x_\tau(t,\varphi))|_{\tau=t}=-2V_1(t,\varphi).$$
From the inequalities $\dot{V}_1(t,\varphi)\leqslant -4V_1(t,\varphi)$ and $\dot{V}_2(t,\varphi)\leqslant -2V_1(t,\varphi)$ it follows that the comparison system for system $(4.1)$ is the system
$$\left\{\begin{array}{l}
\dot{\xi}_1=-4\xi_1\\
\dot{\xi}_2=-2\xi_1\\
\end{array}\right.,\quad \xi=(\xi_1,\xi_2)\in\mathbb{R}^2.\eqno(4.2)$$
It is easy to verify that for the vector function $f(t,\xi)=(-4\xi_1,-2\xi_1)$ the condition $f(t,\xi)\in W$ is satisfied. Let us show that the solutions of system $(4.2)$ are ultimate-uniformly $\mu$-bounded with $\mu_0$-control, where $\mu=(\xi_1)$ and $\mu_0=((\xi_0)_1)$. Consider for this system the Lyapunov function $W(t,\xi)=\xi_1^2$. It is clear that
$b(\|\mu\|)\leqslant  W(t,\xi)\leqslant a(\|\mu\|)$, where $b(r)=a(r)=r^2$, and $\dot{W}_{(4.2)}(t,\xi)\leqslant -c(\|\mu\|)$, where the continuous function $c(r)$, $r\geqslant 0$, satisfying the condition $c(r)>0$ for $r>0$, is defined by the equality $c(r)=8r^2$. From this, by the theorem on ultimate-uniform partial boundedness of solutions with partial control of initial conditions \cite{Lap} (see also \cite{Lap3}), it follows that the solutions of system $(3.5)$ are ultimate-uniformly $\mu$-bounded with $\mu_0$-control. Thus, for
system $(4.1)$ all conditions of Theorem 4.1 for $k=p=1$ are satisfied and, consequently, the solutions of system $(4.1)$
are ultimate-uniformly $(x_1)$-bounded with $\varphi_0^{(1)}$-control.

\vspace{2cm}

\noindent K.S. Lapin - Mordovian State Pedagogical University named after M.E. Evseviev,
Saransk, Russia, e-mail: klapin@mail.ru
\vspace{0.5cm}

\noindent S.V. Lapin - Serov Str., Saransk, Russia,
e-mail: slapin@mail.ru


\begin{thebibliography}{99}

\bibitem{Io2} Yoshizawa T. Liapunovs function and boundedness of solutions. //
Funkcialaj Ekvacioj, 1959, vol. 2, pp. 95-142.
\bibitem{L} Lyapunov A.M. The General Problem of the Stability of Motion (in russian), Kharkov, 1892.
\bibitem{A} Abdullin R.Z., Anapolsky L.Yu., Voronov A.A.,
Zemlyakov A.S., Kozlov R.I., Malikov A.I., Matrosov V.M. The Method
of Vector Lyapunov Functions in Stability Theory (in russian), Moscow: Nauka, 1987.
\bibitem{M} Matrosov V.M. Method of Lyapunov Vector Functions: Analysis of Dynamical Properties of Nonlinear Systems (in russian), Moscow:
Fizmatlit, 2001.
\bibitem{Lap1} Lapin K.S. Partial uniform boundedness of solutions of systems of differential equations with partly controlled
initial conditions. // Differ. Equat., 2014, vol. 50, no. 3, pp. 305–311. DOI: 10.1134/S0012266114030045.
\bibitem{Lap2} Lapin K.S. Lyapunov vector functions and partial boundedness of solutions with partially controlled initial conditions. //
Differ. Equat., 2016, vol. 52, no. 5, pp. 549–556. DOI: 10.1134/S0012266116050037.
\bibitem{Lap} Lapin K.S. Ultimate boundedness with respect to part of the variables of solutions of systems of differential
equations with partly controlled initial conditions. // Differ. Equat., 2013, vol. 49, no. 10, pp. 1246–1251. DOI: 10.
1134/S0012266113100054.
\bibitem{Lap3} Lapin K.S. Poisson Boundedness of
Solutions of Systems of Differential Equations (in russian), Saransk: RITS MGPU, 2022.
\bibitem{Io1} Yoshizawa T. Stability Theory by Liapunov's Second Method. // Publications of the Mathematical
Society of Japan 9, Math. Soc. of Japan, Tokyo, 1966, 223 pp.
\bibitem{Kra} Krasovsky N.N. Some Problems of the Theory of Stability of Motion (in russian), Moscow: Gostekhizdat, 1959.
\bibitem{H} Hale J. Theory of Functional Differential Equations, New York: Springer, 1977.
\bibitem{R} Rouche N., Habets P., Laloy M. Stability Theory by Liapunov’s Direct Method, New York: Springer-Verlag,
1977.
\bibitem{W} Wazewski T. Systemes des equations et des inegalites differentielles ordinaires aux deuxiemes membres monotones et
leurs applications. // Ann. Soc. Polonaise Math., 1950, vol. 23, pp. 112-166.
\bibitem{Lap4} Lapin K.S. Lyapunov–Krasovsky functionals and uniform-ultimate partial
boundedness of solutions with partial control of prehistories. // Differ. Equat., 2026, vol. 62, no. 6, pp. 723-736.
DOI: 10.7868/S3034503026060018.
\bibitem{Rum} Rumyantsev V.V., Oziraner A.S. Stability and stabilization of motion with respect to a part
of variables (in russian), Moscow: Nauka, 1987.
\bibitem{Lap5} Lapin K.S. Higher-order derivatives of Lyapunov functions and partial boundedness of solutions with partially controllable initial
conditions. // Math. Notes, 2017, vol. 101, no. 6, pp. 1000–1008. DOI: 10.1134/S0001434617050273.
\bibitem{Lap6} Lapin K.S. Higher-order derivatives of Lyapunov functions and ultimate boundedness in the sense of Poisson of solutions 
to systems of differential equations. // Siberian Mathematical Journal, 2018, Volume 59, Issue 6, pp. 1100–1104. 
DOI: 10.1134/S0037446618060137.

\end{thebibliography}
\end{document}